\documentclass[aip,graphicx,preprint,author-numerical]{revtex4-1}
\usepackage{TFG}
\usepackage{amssymb}
\usepackage{amsfonts}
\usepackage{amsmath}
\usepackage{graphicx}
\usepackage{dcolumn}
\usepackage{bm}
\usepackage{hyperref}
\usepackage[utf8]{inputenc}
\usepackage[T1]{fontenc}
\usepackage{mathptmx}
\usepackage{xcolor}
\usepackage[normalem]{ulem}

\usepackage[shortlabels]{enumitem}

\graphicspath{{RDS_4_DA-Plots/}}
\hypersetup{
        colorlinks=true,
        linkcolor=cyan,
        citecolor=cyan,
        filecolor=cyan,
        urlcolor=blue,
}
\definecolor{rred}{rgb}{0.55, 0.0, 0.0}
\newcommand{\mg}{\color{blue}}

\begin{document}

\title[Forecast--assimilation (FA) process with unstable dynamics]{
Asymptotic behavior of the forecast--assimilation process with unstable
dynamics}
\author{Dan Crisan}
\address{Department of Mathematics, Imperial College London, 180 Queen's Gate, London SW7 2AZ, UK}
\email{d.crisan@imperial.ac.uk}
\author{Michael Ghil}
\address{Geosciences Department and Laboratoire de M\'et\'eorologie Dynamique (CNRS and IPSL), \'Ecole Normale Sup\'erieure and PSL University, Paris, France}
\email{ghil@lmd.ipsl.fr}
\address{Atmospheric \& Oceanic Sciences Department, University of California at Los Angeles, Los Angeles, CA, USA}
\email{ghil@atmos.ucla.edu}
\date{\today }

\begin{abstract}
Extensive numerical evidence shows that the assimilation of observations has
a stabilizing effect on unstable dynamics, in numerical weather prediction
and elsewhere. In this paper, we apply mathematically rigorous
methods to showing why this is so. Our stabilization results do not assume 
a full set of observations and we provide examples where it suffices 
to observe the model's unstable degrees of freedom.      


\bigskip
\begin{center}
        {\it This paper is dedicated to the memory of Anna Trevisan \\ and to her contributions to data assimilation}
\end{center}
\bigskip 
\end{abstract}

\maketitle

\hypersetup{    ,urlcolor=black
    ,citecolor=blue
    ,linkcolor=blue
    }

\preprint{AIP/123-QED}

\affiliation{Department of Mathematics, Imperial College London, London, UK.}

\affiliation{Geosciences Department and Laboratoire de M\'et\'eorologie Dynamique (CNRS and IPSL), 
                Ecole Normale Sup\'erieure and PSL University, F-75231 Paris Cedex 05, France}
\affiliation{Department of Atmospheric and Oceanic Sciences, University of
                California,Los Angeles, CA 90095-1565, USA}

\begin{quotation}
V. Bjerknes first described weather prediction as an initial-value problem
in 1904 \cite{Bjerknes.1904}. As J. von Neumann and associates started using
computers to implement this idea immediately after World War II, it quickly
became apparent that the requisite initial data available were incomplete 
\cite{Panofsky.1949, Bengtsson.ea.1981}. The appearance of weather
satellites in the 1960s led further on to the concept of time-continuous
assimilation of remote-sensing data \cite{Charney.ea.1969, Ghil.ea.1979}.
Nowadays, data assimilation (DA) is being applied across all the areas of
the climate sciences and much beyond \cite{Carrassietal2020,
Roisin.Beckers.2011, Ghil.Mal.1991, Kalnay.2003, Gottwald.Reich.2021}.
Three crucial issues are still insufficiently well understood: (i) standard
proofs for the convergence of the DA process rely on the stability of the
model dynamics, even in the linear case, while atmospheric and oceanic
dynamics are famously unstable and chaotic \cite%
{Ghil.Chil.1987,Kalnay.2003,Lorenz.1963a}; (ii) data availability over time
appears to successfully compensate for insufficient instantaneous coverage
in space \cite{Bengtsson.ea.1981,Ghil.Mal.1991, Titi.ea.DA.2019}; and, last
but not least, (iii) it appears that observations of a model's unstable
manifold are sufficient for the convergence of the time-continuos
forecast--assimilation cycle \cite{Carrassi.ea.2008, Trevisan.Uboldi.2004}.
The present paper uses concepts and methods from  
the stochastic calculus \cite{KaratzasShreve}, random dynamical systems \cite{Arnold.1998,Caraballo.Han.2017,Crauel.Flandoli.1994}
and nonlinear filtering \cite{bc, CR}
to achieve several significant steps in
clarifying all three of these issues. 
\end{quotation}

\section{Introduction and motivation}

\label{sec:intro}

\subsection{The forecast--assimilation cycle in meteorology}

\label{ssec:background}

A key metaproblem of data assimilation (DA) in atmospheric, oceanic and
climate problems is to show that sequential filters of various degrees of
sophistication are stable, and that they converge to solutions with suitable
properties that lie sufficiently close to the observations, such as they
are, in the case of the \textit{unstable dynamics} that characterizes these
problems. Moreover, a full solution to this metaproblem should allow one to
compare, with relative ease, the efficiency and accuracy of several filters.

Heuristically, the motivation for this metaproblem being soluble is the
success of practical DA methods in numerical weather prediction (NWP), and
in related oceanographic and climate problems, in keeping track of a
system's observed state \cite{Ghil.Mal.1991, Ghil.Todling.1996}. We outline
herein some simple ideas of why sequential filters do have a chance of being
stable and convergent, even in the presence of dynamically unstable modes.

As the importance of DA methodology in the climate sciences and an 
increasing number of other areas, all the way to finance, has been growing rapidly, 
DA has attracted more and more attention in relevant areas of the mathematical 
sciences \cite{Asch.Bocquet.2016, bc, rc, Stuart.ea.2015}. As we shall argue 
further below, this increase of attention has not covered as yet the issues of
instability of the basic dynamics that one wishes to track nor that of partial
observations. The main point of this paper is to take substantial steps in
addressing these two issues with the same degree of rigor as that used so
far in addressing the considerably simpler situation of stable dynamics
and complete measurements.

To lead up to the full complexity of the setting involved, consider first, for simplicity, 
a scalar, linear DA problem in continuous time, using ``sloppy'' notation. 
The model equation is 
\begin{equation}
\dot x = Ax+u(t),  \label{model}
\end{equation}%
and the observation equation is 
\begin{equation}
\dot y = Hx+v(t),  \label{obs}
\end{equation}%
where $(u,v)$ are noises that are white in time $t$ and have ``nice'' densities
--- i.e., centered and with finite variance, e.g., Gaussian --- while the
dot notation $(x, y)^\cdot$ stands for the time derivative. We keep caps $%
(A,H)$ for the model and observation operators because the scalar case is
supposed to be just shorthand for the vector--matrix case. The fully 
nonlinear and high-dimensional problems arising in the actual applications 
involve also much more complex noise processes, of course; see,
for instance, the work of C. Nicolis and coworkers
in ref. \cite{Nicolis.ea.2009} and further references therein.


The forecast--assimilation (FA) process for $\hat{x}$, the best linear
unbiased estimate of $x$, obeys\cite{Kalman.1960, Kalman.Bucy.1961} {%
\begin{equation}  \label{FA}
\dot{\hat{x}} = A\hat{x} + K(y - H\hat{x}).
\end{equation}
Here $K$ is a weight matrix, which equals the Kalman-Bucy optimum \cite%
{Gelb.1974, Jazwinski.1970} in the Gaussian-noise case, and $y - H\hat{x}$
is the {\em innovation vector} that equals the difference between the actually
observed value of $x$ and the one expected by forecasting this value based
on past observations.

 The intuitive
motivation to hope for convergence of this FA cycle to the true evolution ---
or for its synchronization with the observations\cite{Duane.ea.2017} --- is
simply re-writing equation (\ref{FA}) as 
\begin{equation}
\dot{\hat{x}}=(A-KH)\hat{x}+Ky.  \label{sync}
\end{equation}%
This equation exhibits the new, and hopefully stabler, dynamics $(A-KH)$ of
the FA process vs. the original, pure-evolution dynamics $A$. It also suggests
using random dynamical system (RDS) theory\cite{Arnold.1998,
Crauel.Flandoli.1994} for the FA problem, given the time-dependent forcing
by $Ky$, in which the observations are subject to random errors.
An exhaustive presentation of RDS theory is given in the Ludwig
Arnold monograph \cite{Arnold.1998}, which makes, however, somewhat 
difficult reading for the non-specialist. More accessible presentations 
for DA practitioners can be found in refs. \cite{Caraballo.Han.2017,Charo.ea.2021,GCS.2008}. 

One has to show that the nonlinear, multidimensional --- and possibly even
infinite-dimensional FA evolution, as in the generalization of Eq.~%
\eqref{sync} to partial differential equations (PDEs) --- is stable \cite%
{Carrassi.ea.2008, Trevisan.ea.2010} even in the presence of dynamic
instabilities. 
For the way that it might still suffice, in the presence of dynamic
instabilities, to have dim$\{y\}\le $ dim$\{x\}$, see numerical results
for simplified atmospheric and oceanic models in ref. \cite{Ghil.Todling.1996} 
and in Section 2 of the Ghil (1997)\cite{Ghil.1997} review paper, for instance. 
Frank and Zhuk \cite{Frank.Zhuk.2018} did obtain
a mathematically rigorous result along these lines for a deterministic
system of nonlinear ordinary differential equations.

We restrict ourselves here to the finite-dimensional case: in the
operational practice of numerical weather prediction (NWP), the PDEs
governing atmospheric and oceanic flows are discretized in physical space
--- using finite differences, finite elements or spectral and
pseudo-spectral methods \cite{Kalnay.2003,Roisin.Beckers.2011}. These days,
the number of resulting finite-difference equations in time is very large
--- up to order of $10^8$--$10^9$ --- but still finite and will stay so for
the foreseeable future.

The nonlinearity of atmospheric and ocean dynamics \cite{Ghil.Chil.1987,
Ghil.Luc.2020, Lorenz.Book.1967} compels us to deal not just with the mean
and variance of the estimated state $\hat x$, as in the classical
Kalman-Bucy filter \cite{Kalman.1960, Kalman.Bucy.1961}, but with the 
entire probability distribution function (pdf) of the state $x$, conditioned on the observations $z$.
This pdf may be multimodal or, more generally, non-Gaussian \cite{Asch.Bocquet.2016, CR} 
and include the presence of long tails due to extreme events \cite%
{Chavez.ea.2015, Ghil.ea.ExEv.2011}. 

To fully describe the pdf of the observed state $x(t)$ given the data $z(t)$, one needs more than this pdf's mean and variance. Indeed, absent the linearity and Gaussianity assumptions, the system of equations satisfied by the mean and the variance of the pdf of the observed state is no longer closed and we need more that just these two quantities to describe the FA process and compute its evolution. This point of view is, by now, widely shared by the operational NWP community \cite{Bocquet.ea.2010, Carrassietal2020, Leeuwen.2009}.


In this paper, we assume that the model $x$ is a, possibly unstable, stochastic process
explicitly defined in section \ref{sec:framework}. We also assume that
the model is observed only partially and that the observations arrive
continuously in time. The latter assumption is consistent with the already
mentioned continuous flow of observations in the satellite era \cite%
{Charney.ea.1969, Ghil.ea.1979}.

In fact, shortly after the advent of meteorogical satellites in the late
1960s, Charney et al. \cite{Charney.ea.1969} formulated the conjecture that
a knowledge of the continuous time history of the atmospheric temperature
field will allow one to determine the other state variables, in particular
the winds. M. Ghil and coauthors provided analytical arguments for the
correctness of this ``Charney conjecture'' in two-dimensional
(2-D) geophysical fluid dynamics (GFD) models \cite{Ghil.ea.1977} and
documented its usefulness numerically with time-continous DA of actual
remotely sensed temperatures in a fairly realistic NWP-type,
three-dimensional (3-D) atmospheric model \cite{Ghil.ea.1979}. E. Titi and
coauthors provided rigorous proofs in a purely deterministic setting for
both 2-D and 3-D models of GFD interest \cite{Farhat.ea.2015,
Titi.ea.DA.2019}.


Figure~\ref{fig:Ghil1981} here illustrates the FA process's evolution in the 
presence of observations that are partial in both their nature --- i.e., 
temperature vs. winds --- and their spatial coverage --- ocean vs. land 
in this simple example. The figure represents DA results using the Kalman-Bucy
filter\cite{Gelb.1974,Jazwinski.1970,Kalman.1960,Kalman.Bucy.1961}
(hereafter KF) for a linear, mid-latitude shallow-water model in one space
dimension with a simplified geometry \cite{Ghil.ea.1981}. In this geometry,
there are two data-rich regions of equal size that stand for the North
American and Eurasian land masses, alternating with two data-poor regions of
the same size that stand for the North Atlantic and North Pacific.

In the particular numerical experiment selected here for illustration purposes, 
all three model variables --- namely the geopotential height $\phi = gh$ of the
free surface, where $h$ is the actual height and $g$ the acceleration of
gravity, along with the cartesian velocity components $(u, v)$ --- were
available at the so-called synoptic times of noon and midnight GMT over the
land areas, while no data at all were available over the ocean areas. In
this case, all three curves for $u (t), v(t)$ and $\phi (t)$ --- as well as
for the energy $E = u^2 + v^2 + \phi^2/\Phi$, where $\Phi$ is the
equilibrium value of $\phi$ about which the equations are linearized ---
have exactly the same behavior.

Over the data-rich land, the error drops sharply at the first observing time,
12~h after the start, and it grows parabolically between each observing time
and the next one, due to the KF's Gaussian and white-in-time model noise. 
The parabolic growth of the root-mean-square (RMS) error in this case
is due to the additive model noise, as in the case of scalar Brownian motion 
\cite{Wax.1954}, in which it is the variance of the process that is
proportional to time $t$; compare with the dashed red line in Fig.~\ref{fig:FA_Cycle}
below. Over land, though, the evolution of the RMS error
asymptotes very quickly, in 1--2 model days, to one in which the FA error
level oscillates around that of the observational noise, which equals
roughly 0.9 in nondimensional units.

\begin{figure}[!ht]
\centering
\includegraphics[width=0.65\linewidth]{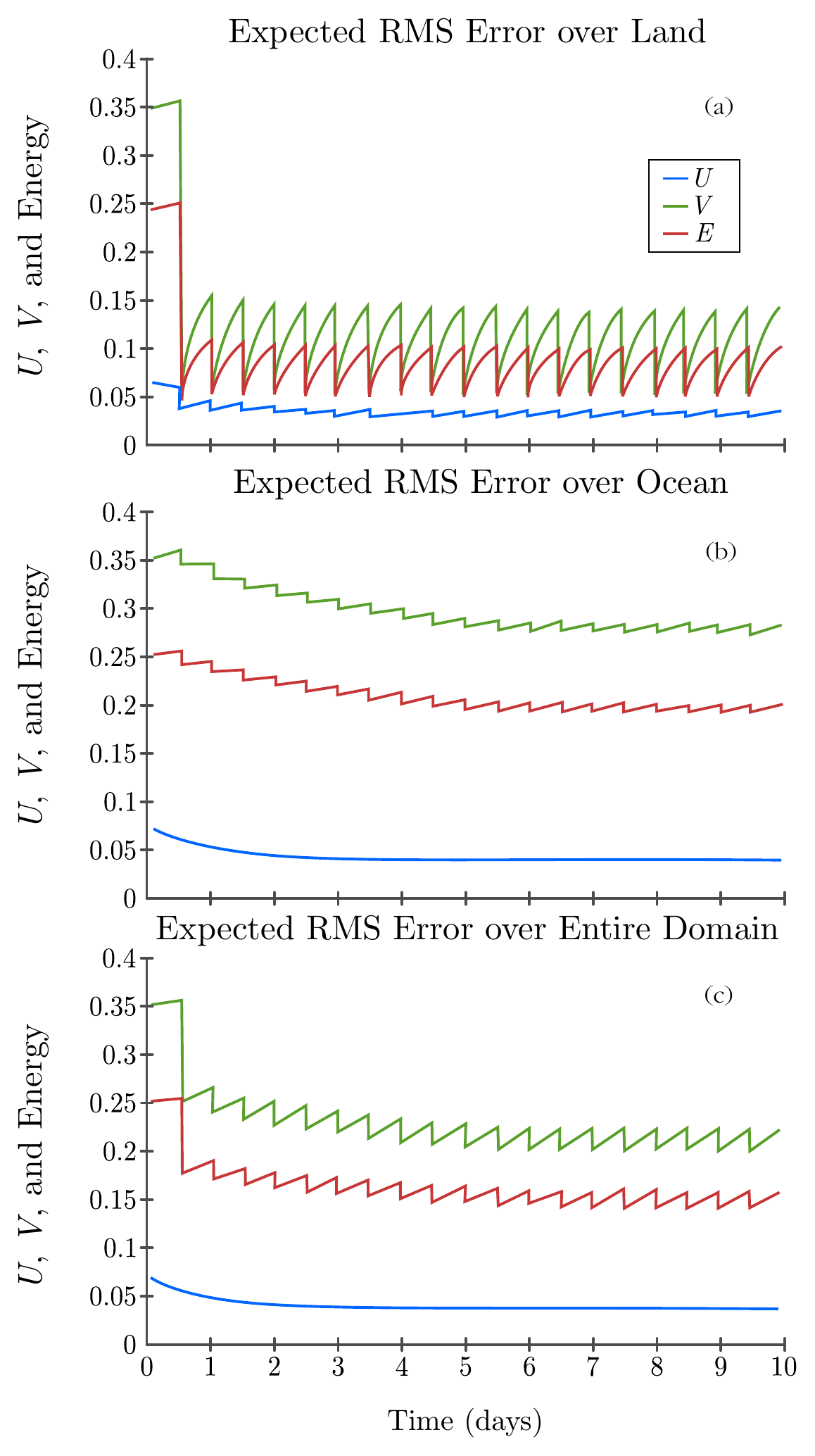}
\caption{Typical results of a forecast--assimilation (FA) cycle. 
        (a) Expected root-mean-square (RMS) error over land; 
        (b) expected RMS error over the ocean; and 
        (c) expected RMS error over the entire model domain. 
        This figure is based on Fig.~6 in Ghil et al. \protect\cite{Ghil.ea.1981}, 
        with permission, and E. Bach provided the version herein.}
\label{fig:Ghil1981}
\end{figure}

Over the data-poor oceans, the RMS error still decreases, due to the
advection of information from the land areas by the mean westerly winds, $U
= - \partial \Phi/\partial y$. But this error decrease is slower than over
land, the asymptotically periodic behavior is only reached after 4--5
model days, and the mean values of RMS errors stay above the observational
noise. Finally, the RMS error behavior over the entire area is essentially a
weighted average of the results over land and over ocean.

This type of behavior was modified by the advent of time-continuous
satellite data over the oceans, as shown, for instance, by the work of Halem
et al.\cite{Halem.ea.1982}, as part of the DA studies associated with the
Global Atmospheric Research Experiment (GARP). While the Ghil et al. \cite%
{Ghil.ea.1981} model was one-dimensional, linear, stable and had only a
rather small spatial resolution, that of Halem et al.\cite{Halem.ea.1982}
and many others were fully 3-D, nonlinear, unstable, and had rather high
resolution by the standards prevailing at that point in time.

Figure~5 in Halem et al.\cite{Halem.ea.1982} (not shown here) clearly
indicates the improvement in 6-hr forecasts over the Western U.S. from
initial states that do use the time-continuous satellite data over the North
Pacific vs. those that use only the conventional data available over land. 
The mechanism that advects information by the westerly winds --- or misinformation,
from ocean to land, in the case of the conventional observing network --- is
clearly still working, in spite of the unstable and nonlinear dynamics of
the forecast model and of the ad-hoc sequential-estimation method in ref.
\cite{Halem.ea.1982}. The latter was a successive-correction method 
\cite{Cressman.1959,Ghil.ea.1979} in the Halem et al.\cite{Halem.ea.1982} 3-D
model, while it was a KF in the Ghil et al. \cite{Ghil.ea.1981} 1-D model.

An idealized version of the error components that play a role in the 
FA process are illustrated in Fig.~\ref{fig:FA_Cycle} and discussed further 
in {\mg Appendix~I}. The major point
of the figure and of the associated appendix is that high-end, operational
NWP models have various sources of instability, and that NWP would be
impossible if DA did not stabilize the FA dynamics and thus achieve
real-time low-error forecasts. The main purpose of the present paper is to
justify rigorously the numerically observed fact that this is, indeed, the
case.

\begin{figure}[!hb]
\centering
\includegraphics[width=0.7\linewidth]{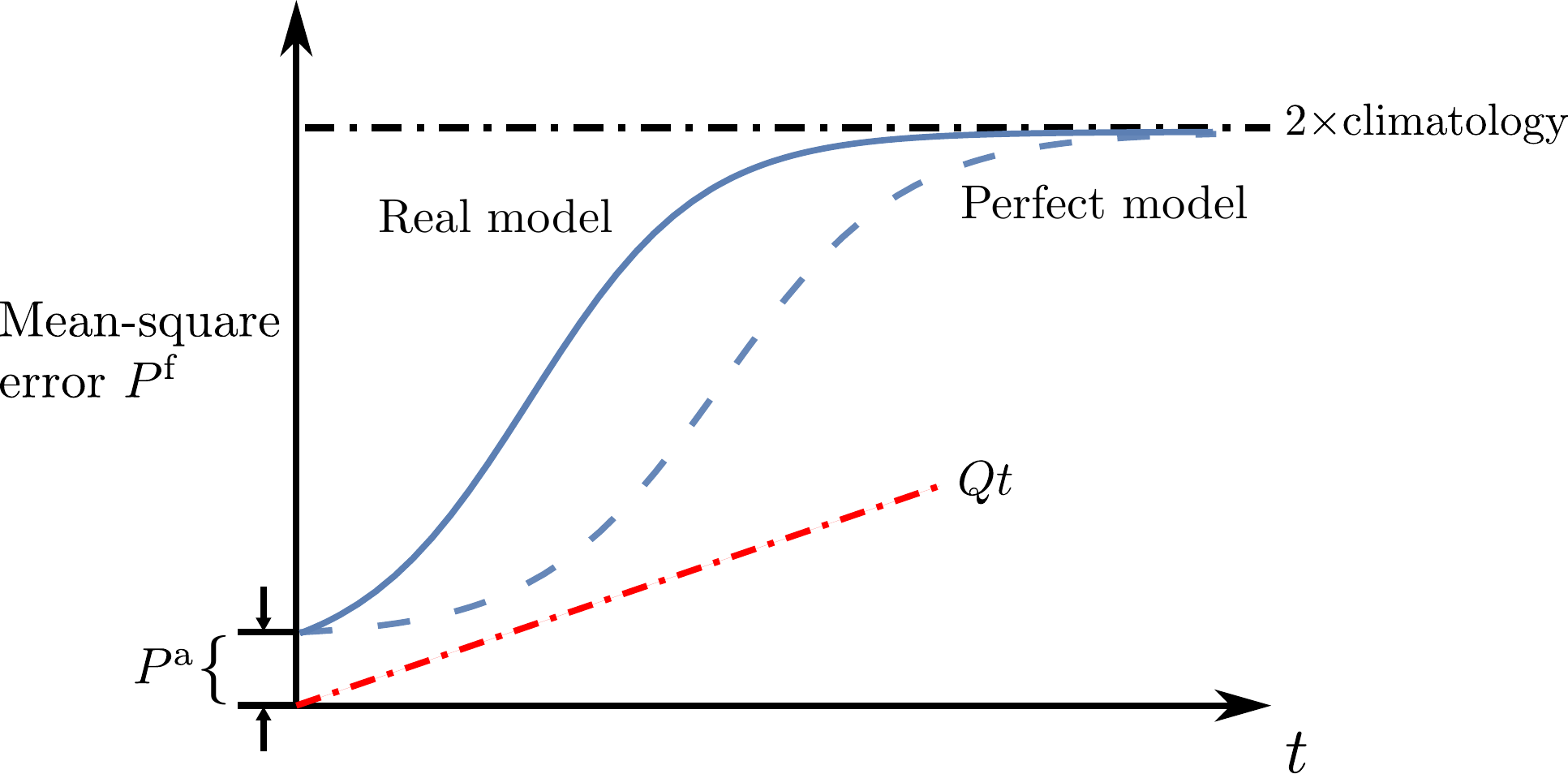}
\caption{Schematic diagram of forecast error growth in a numerical weather
prediction (NWP) model; see Appendix~I for explanations. The perfect-model
(dashed blue) and real-model (solid blue) curves here should be compared
with the blue curve in the two panels of Fig.~\protect\ref{fig:error_models}
there, for $S = 0$ and $S \neq 0$, respectively. The straight red line
(dash-dotted) labeled $Qt$ here represents the linear growth of forecast
error variance due to additive white noise, as is the case in Fig.~\ref{fig:Ghil1981}(a),
between updates and over land,
while the Lorenz \cite{Lorenz.1982} and Dalcher and Kalnay \cite%
{Dalcher.Kalnay.1987} models in the appendix only assume a constant
deterministic model error $S$. E. Bach kindly provided this figure.}
\label{fig:FA_Cycle}
\end{figure}

\subsection{A rough sketch of the mathematical formulation}

\label{ssec:formulation}

In the present paper's context, the ``best estimate'' of the model $x$ is
its conditional distribution with respect to the observational data $z(t)$
available up to the current time $t $. We will denote this probability
distribution by $\pi (t) \equiv \pi_t$ and the distribution of the model $x$
in the absence of any observational data by $p(t) \equiv p_t$. The pdf $p$
is called the prior distribution and $\pi$ the posterior distribution, where
we dropped for simplicity the dependence on time $t$. The two distributions 
$\pi$ and $p$ can be viewed as dynamical systems that both evolve in the
infinite-dimensional space of probability measures $\mathcal{P}$ over 
the model state space $\mathbb{R}^{d}$. Put in simple words, our results 
refer to the ideal, truly optimal filter and not to any specific approximations 
thereof, such as the extended Kalman filter (EKF~\cite{Ghil.Mal.1991}) or 
ensemble Kalman filter (EnKF~\cite{Bocquet.ea.2010, Carrassietal2020, 
Leeuwen.2009}) or any of the variational methods~\cite{Ghil.Mal.1991,Bocquet.ea.2010}. 

For the comfort of the interested reader coming from the NWP community 
or from other areas where DA is being used --- or its use is being contemplated 
--- we are summarizing in Table {\color{blue} I}  the correspondence of 
key symbols and terms in this paper vs. their counterparts in NWP use.
\begin{table}[h] \label{tab:notation}
        \begin{tabular}{|l|l|l|}                
                \hline
                
                Notation &Data Assimilation Language &  Stochastic Filtering Language \\             
                \hline
                
                $p_t$ &Probability distribution of pure forecast & Prior distribution of the signal \\           
                \hline              
                $\pi_t$ & Probability distribution of the analysis &  Conditional distribution of the signal \\                   
                \hline
                
                $ \lambda_t^f$ & Deterministic forecast model $M$ \cite{Ide.ea.1997} &  Push-forward operator \\
                \hline        
        \end{tabular}
        \caption{Correspondence of symbols and terms in this paper vs. their counterparts in NWP}
\end{table} 
Depending on the model for $x \in \mathbb{R}^{d}$, the prior system may be
unstable: Starting the prior system $p$ from two different initial conditions, the $p_0$  
that results in $p_t$ and $\mu \neq p_0$, and using the 
same forecast operator for both
produces two probability measures, $p_t$ and $p^\mu_{t}$, that will diverge from 
each other in time, in the sense that the Wasserstein distance $W_2$ between the two
tends to infinity. The distance $W_2$ between two probability measures 
$\mu$ and $\nu$ is given by 
        \begin{equation}\label{eq:W2}
                {W_{2}(\mu ,\nu ):=\left( \inf_{\gamma \in \Gamma (\mu ,\nu )}\int_{{{%
                                                \mathbb{R}}^{d}}\times {{\mathbb{R}}^{d}}}|x-y|^{2}\,\mathrm{d}\gamma
                        (x,y)\right) ^{1/2},}
        \end{equation}%
        where ${\Gamma (\mu ,\nu )}$ denotes the collection of all measures on ${{%
                        \mathbb{R}}^{d}}\times {{\mathbb{R}}^{d}}$ with marginals $\mu $ and $\nu $ on the first and second factor, respectively. 

This distance is often called in applications the 
“earth mover's distance,” following the original motivation of G. Monge\cite{Monge.1781}.
Ghil \cite{Ghil.2015} originally proposed the idea of using the Wasserstein distance
in the context of the climate sciences as a way to generalize the traditional concept 
of equilibrium climate sensitivity \cite{Ghil.Luc.2020} in the presence of 
a time-dependent forcing, such as seasonal or anthropogenic forcing.
Further details are given in Appendix II.A and considerably more information on the definitions and methods used herein can be found in the Panaretos and Zemel monograph \cite{Panaretos.Zemel.2020}.


As stated already several times, we concentrate on the prior process's $p_t$ being unstable,
i.e., starting it from the initial condition $\mu \neq p_o$ will lead to divergence of the 
trajectory $p^\mu_{t}$ from $p_t$ in Wasserstein distance $W_2$. We show in 
Sec.~\ref{sec:wd} that, to the contrary, starting the posterior process $\pi^\mu_{t}$ 
from the initial condition $\mu \neq p_o$ and evolving it with the same FA operator 
as that used to evolve $\pi_{t}$ will generate a probability measure $\pi^\mu_{t}$ that 
will keep the  $W_2$ distance between $\pi^\mu_{t}$ and $\pi _{t}$ bounded in
expectation. Moreover, in the linear case, we show that  the $W_2$ 
distance between $\pi^\mu_{t}$ and $\pi _{t}$ actually tends to zero.

Mathematically, we consider the forecast operator $\lambda _{t}^{f }$ and
the FA operator $\lambda _{t}$ associated with the FA process, respectively,
which are defined by  
\begin{subequations}
\label{eq:processes}
\begin{align}
& \lambda _{t}^{f }: \mathcal{P}(\mathbb{R}^{d})\rightarrow \mathcal{P}( 
\mathbb{R}^{d}),\ p_{t}=\lambda _{t}^{f }p_{0},  \label{eq:fcst} \\
& \lambda _{t}: \Omega \times \mathcal{P}(\mathbb{R}^{d})\rightarrow 
\mathcal{P}(\mathbb{R}^{d}),~\pi _{t}\left( \omega \right) =\lambda
_{t}\left( \omega \right) \pi _{0},  \label{eq:FA}
\end{align}
and we show that, under certain conditions, 
\end{subequations}
\begin{equation}
\sup_{t\in \lbrack 0,\infty )}\mathbb{E}\left[ W_2\left( \lambda
_{t}\mu,\lambda _{t}\pi_{0}\right) \right] <\infty ,  \label{eq:expl}
\end{equation}%
whilst possibly $\lim_{t\rightarrow \infty }W_2( \lambda _{t}^f\mu,\lambda _{t}^fp_{0}) =\infty$. 
Moreover, in the linear case, we will show that $\lim_{t\rightarrow \infty
}W_2\left( \lambda _{t}\mu,\lambda _{t}\pi_{0}\right) =0$. 

Note that, in Eqs.~\eqref{eq:processes} above, both the push-forward operators $\lambda _{t}^{f }$ and $\lambda _{t}$ are stochastic processes and
act on the full probability measures $p$ and $\pi$, respectively. In the present  framework, we do
not limit ourselves just to the mean and variance of the state $x(t)$, as is
the case for the linear KF \cite{Kalman.1960, Kalman.Bucy.1961}.  The distance $W_2\left( \lambda
_{t}\mu,\lambda _{t}\pi_{0}\right)$ between the RDS $\lambda
_{t}\mu$ starting from $\mu$ and the RDS $\lambda
_{t}\pi_0$ starting from $\pi_0$  is of course random: It depends on the realization of the observation process, 
and hence the expectation in Eq.~\eqref{eq:expl} is taken with respect to the pdf of  this process.


We will prove rigorously, in the precise sense described above, that the incorporation
of observational data into the FA process does indeed have a stabilizing
effect on unstable dynamics, as shown abundantly by NWP practice and
suggested, in particular, by the work of Anna Trevisan and her collaborators
on assimilation in the unstable manifold (AUS)\cite{Carrassi.ea.2008,
Trevisan.Uboldi.2004, Trevisan.ea.2010}. The latter work was one source of
inspiration for the present paper.

Another source was RDS theory, according to the arguments presented by one
of the authors (MG) at the ``Symposium Honoring the Legacy of Anna
Trevisan'' held in Bologna, Italy in October 2017. A more specific source of
inspiration for the rigorous mathematics herein was the paper of L. Arnold 
\cite{Arnold.1990} on stabilization by noise, which was presented to MG by
Franco Flandoli during the trimester on ``The Mathematics of Climate and the
Environment'', held at the Institut Henri Poincar\'e %
in Paris in Fall 2019. 

It turned out, however, rather quickly that the
latter paper's arguments could not be applied directly to the DA problem at
hand, since the deterministic component of the process under study there is
linear, and we did not see how to extend Arnold's arguments \cite
{Arnold.1990} to fully nonlinear processes. 
 Moreover, the latter arguments are only valid for finite-dimensional
systems, while we are dealing here with RDSs that evolve in the infinite-dimensional space of probability
measures $\mathcal{P}(\mathbb{R}^d)$.

The layout of the paper is dictated by the intent to bring the two
communities --- of DA practitioners, on the one hand \cite%
{Asch.Bocquet.2016, Bengtsson.ea.1981, Carrassi.ea.2008, Ghil.Mal.1991} ---
and of the rapidly increasing numbers of applied mathematicians and
physicists interested in DA on the other \cite{CR, Leeuwen.ea.2015,
Stuart.ea.2015} --- closer together. Hence, after this fairly long
introduction, we describe in Sec.~\ref{sec:framework} the precise
mathematical framework that is used herein. The main rigorous results are
outlined in Sec.~\ref{sec:wd} and further details on definitions and proofs
appear in {\mg Appendix~II.} Conclusions and some thoughts on further work appear
in Sec.~\ref{sec:concl}.

\section{Mathematical framework} \label{sec:framework}

\subsection{Prior results}  \label{ssec:prior}

The two stabilization theorems presented here, Theorems~\ref{nonlinear} and \ref{linear}, are related to 
and do use in their proofs certain arguments from the existing 
stability results in the nonlinear filtering literature, e.g., refs. \cite{Budhiraja2011, Moral.1996, dmlm, uc2, uc3, kun, op, picard1991,Handel.2009}. Such results have not, by and large, matched the objectives of the applied DA community when studying the asymptotic behavior of the FA process. 

To be more precise, the applied DA community is interested in results for a forecast cycle that is unstable --- as is the 
case in meteorology and oceanography --- and for which applying DA has the mysterious but 
salutary effect of stabilizing the FA process. In addition, the DA community has to rely, 
typically, on results where only a small subset of the forecast cycle's degrees of freedom can be 
observed. Without being comprehensive, of course, we give here a classification of the conditions under which the existing results on filter stability for \emph{nonlinear} forecast processes hold.       

\begin{itemize}
        
        \item The forecast process is assumed to be ergodic or to have good mixing properties, e.g. Atar\cite{atar2011}, Atar and Zeitouni\cite{az1997}, Budhiraja\cite{Budhiraja2011}, Chigansky and Liptser\cite{cl2004},   Chigansky, Liptser and Van Handel\cite{clh2011}, Del Moral and Miclo\cite{dmlm},  and Del Moral, Doucet and Singh\cite{uc3}. Please note that ergodicity in these papers  is meant in the sense of convergence of the pdf of the push-forward process $\lambda^f$ as defined  in Eq.~\eqref{eq:fcst}, not in the pathwise sense of RDSs \cite{Arnold.1998, Checkrounetal2011, GCS.2008}. For example, the stochastic Lorenz 
model\cite{Checkrounetal2011} is ergodic in the pathwise sense of RDSs, but does not have a stable push-forward process $\lambda^f$. Note that ergodicity of the pdf of the push-forward process $\lambda^f$ does imply its stability. For further results on ergodicity of the push-forward process, see, for instance, Bakry, Cattiaux and Guillin \cite{Bakry.ea.2008}.   
        
        \item The forecast cycle is fully observed, e.g., Picard\cite{picard1991}, 
        Crisan and Heine\cite{ch2008}, Stannat\cite{stannat2011}, 
        Chigansky, Liptser and Van Handel\cite{clh2011}, and Van Handel\cite{Handel.2009}.
        
        \item The forecast process has a drift term which is the gradient of a convex function, or a perturbation thereof,  and  its noise term is strictly elliptic, e.g., Stannat\cite{stannat2005, stannat2011}. This is a very restrictive class of forecast processes, which  are expected to be  stable; see Bakry \& all\cite{bb2008}. Specific drifts belonging to this class are used in the classical Metropolis-adjusted Langevin algorithm known to converge faster to the diffusion's invariant measure; e.g., Roberts and Tweedie\cite{rt1996}.      
        
\end{itemize}

To try to answer the question raised by the applied community, we no longer insist on proving 
that the FA process is (exponentially) stable in the sense advocated by the theoretical 
community. We relax the definition and only require that the FA process initialised from the 
wrong distribution does not diverge too strongly from the correctly initialized FA process, even 
when the forecast process does so; this is why the ergodicity assumption for the forecast process is 
not useful. 

We substantially strengthen, however, the stabilization result, in the sense that we 
want to control the mean and the second moment of the FA process. Again this is needed for 
practical reasons. The practitioners want to know that the pointwise estimate of their algorithm 
of choice does not diverge from the theoretical mean of the FA process. But they also want to 
know that their error bars are not too different from the theoretical ones.            

Furthermore, the results presented herein have the advantage that the signal $X(t)$ is not required to be fully observable. In particular, the dimension $n$ of the observation $Y(t)$ and the dimension $d$ of the signal $X(t)$, in the notation of Sec.~\ref{sec:framework}, do not need to coincide; in operational NWP and many other applications, $d \ll n$ \cite{Asch.Bocquet.2016,Bengtsson.ea.1981}. 

Heuristically we need to be able to observe all the "unstable" directions, as suggested by A. Trevisan and coworkers
\cite{Carrassi.ea.2008, Trevisan.Uboldi.2004, Trevisan.ea.2010}, who gave several fairly realistic examples of this idea working quite well. The connection with our results is provided by some simple illustrative examples in  {\mg  Appendix II.D}. 


In the linear case, the stability of the FA process is better understood;   
see, for instance Ocone and Pardoux\cite{op},  Picard\cite{picard1991}, Van Handel\cite{Handel.2009}, and
Stannat\cite{stannat2005}.  
The stability of the associated Riccati equation has been studied in Bishop and Del Moral\cite{bmm1,bmm2}. The stabilization result obtained in Theorem \ref{linear} below for the linear signal is weaker than many of the  existing results discussed above. For example, no (exponential) rates of convergence are deduced herein. This situation clearly leaves considerable room for proving stronger results, given the hypotheses of Theorem \ref{linear}.

\subsection{The present setting}  \label{ssec:setting}
Having presented in the previous subsection a quick review of previous mathematical results
on the stability of the FA process, we proceed now by introducing the setting of our two theorems
in the next section.
Let $X=(X^{i})_{i=1}^{d}$ be the solution of the following stochastic
differential equation driven by a $p$-dimensional Brownian motion process $%
V=(V^{j})_{j=1}^{p}$, 
\begin{equation}
X_{t}=X_{0}+\int_{0}^{t}f\left( X_{s}\right) \mathrm{d}s+\int_{0}^{t}\sigma
\left( X_{s}\right) \mathrm{d}V_{s}.  \label{signal}
\end{equation}%
Here, we assume that $f=(f^{i})_{i=1}^{d}:\mathbb{R}%
^{d}\rightarrow \mathbb{R}^{d}$ and $\sigma =(\sigma ^{ij})_{i=1,\ldots
,d,j=1,\ldots ,p}:\mathbb{R}^{d}\rightarrow \mathbb{R}^{d\times p}$ are
globally Lipschitz. This will ensure that the equation \eqref{signal} has a unique solution.

As stated in Sec.~\ref{ssec:formulation}, we are interested in tracking the evolution of 
the full pdf of the prior and posterior processes, namely $p_t$ and $\pi_t$, respectively.
To do so, we recall that the process $X$ is a diffusion process with infinitesimal generator given by 
\begin{equation*}
A\varphi =\sum_{i,j}\frac{1}{2}a_{ij}\partial _{i}\partial _{j}\varphi
+\sum_{i}f_{i}\partial _{i}\varphi ,
\end{equation*}%
where $a_{ij}=\sum_{k}\sigma _{ik}\sigma _{jk}$. The prior distribution of $X_{t}$ 
is also called its law in the context of filtering \cite{CR}. 

For an arbitrary measurable function $\varphi :\mathbb{R}^{d}\rightarrow \mathbb{R}$ 
that is integrable with respect to the law of $X_{t}$, one has
\begin{equation*}
p_{t}(\varphi )=\mathbb{E}\left[ \varphi \left( X_{t}\right) \right] .
\end{equation*}%
By restricting $\varphi $ further to lie in a suitably chosen space of functions denoted by $%
\mathcal{D}(A)$, 
\begin{equation}
p_{t}(\varphi )=p_{0}(\varphi )+\int_{0}^{t}p_{s}(A\varphi )\mathrm{d}s.
\label{priordistributionequation}
\end{equation}%
Let $\lambda _{t}^{f}:\mathcal{P}(\mathbb{R}^{d})\rightarrow \mathcal{P}(%
\mathbb{R}^{d})$ be the push-forward operator associated with equation %
\eqref{priordistributionequation}. In particular, for arbitrary $\mu \in 
\mathcal{P}(\mathbb{R}^{d})$, the map $t\rightarrow \lambda _{t}^{f}\mu
=:p_{t}^{\mu }$ is the solution of \eqref{priordistributionequation} with
the initial condition $\mu $.

Let $Y$ be an $n$-dimensional process, where $n$ is the number of
observations that are taken to be one-dimensional and have measurement
noises modelled by independent Brownian motions: 
\begin{equation}
Y_{t}^{i}=Y_{0}^{i}+\int_{0}^{t}h^{i}(X_{s})\,\mathrm{d}%
s+W_{t}^{i},~~~i=1,...,n,  \label{observation}
\end{equation}%
and $h_{i}\left( x\right) $ are the corresponding observation operators,
which generalize the observation matrix from the linear framework; see
{\mg Appendix II.C} for details.

Let $\pi =\{\pi _{t},t\geq 0\}$ be probability measure--valued process that
gives us, at time $t\geq 0$, the conditional distribution of the signal $%
X_{t}$ given the observations accumulated up to time $t$, $\{Y_{s},\ s\in
\lbrack 0,t]\}$. It is this process that we call the FA process in our rigorous mathematical 
setting. In other words, $\pi _{t}$ satisfies 
\begin{equation*}
\pi _{t}(\varphi )=\mathbb{E}\left[ \varphi \left( X_{t}\right) |Y_{s},\
s\in \lbrack 0,t]\right] ,
\end{equation*}%
where $\varphi :\mathbb{R}^{d}\rightarrow \mathbb{R}$ is an arbitrary
measurable function that is integrable with respect to the law of $\pi _{t}$. The
FA process satisfies the following stochastic partial
differential equation, formulated here in the integral form 
\begin{eqnarray} 
\pi _{t}(\varphi ) &=& \ \pi _{0}(\varphi )+\int_{0}^{t}\pi _{s}(A\varphi )\,%
\mathrm{d}s+\int_{0}^{t}\left( \pi _{s}(\varphi h^{\top })-\pi _{s}(h^{\top
})\pi _{s}(\varphi )\right) (\mathrm{d}Y_{s}-\pi _{s}(h)\,\mathrm{d}s), \label{ks0} \\
&=& \ \pi _{0}(\varphi )+\int_{0}^{t}\pi _{s}(A\varphi )\,\mathrm{d}%
s+\int_{0}^{t}\left( \pi _{s}(\varphi h^{\top })-\pi _{s}(h^{\top })\pi
_{s}(\varphi )\right) \mathrm{d}I_{s}, \label{ks1} 
\end{eqnarray}%
for any test function $\varphi \in \mathcal{D}(A)$. Here, $I$ is
the innovation process, defined as  
\begin{equation*}
I_{t}^{i}=Y_{t}^{i}-\int_{0}^{t}\pi _{s}\left( h_{i}\right) \mathrm{d}%
s,~~~i=1,...,n,~~~t\geq 0,
\end{equation*}%
which is the rigorous analog of the innovation vector in Eq.~\eqref{FA}
of Sec.~\ref{ssec:background}. The innovation process is a Brownian motion, see e.g., Chapter 3 in Bain\&Crisan\cite{bc}. In particular it is a martingale and stochastic integrals with respect to the  innovation process are easier to manipulate. In particular upper bounds for the stochastic integrals with respect to martingales are easier to obtain than those  for  stochastic integrals with respect to general  semi-martingales. The observation process is a semi-martingale (it is a Brownian motion plus a drift term) so harder to handle. 

More on stochastic partial differential equations and their difficulties can be found, for instance, in Liu and R\"{o}ckner\cite{lr}, or Rozovsky and Lototsky\cite{rl}. 

In order to analyze the asymptotic behavior of the FA process, we can recast 
the solution of  \eqref{ks0} as an RDS. More
precisely, there exists a measurable map $\lambda :[0,\infty )\ \times 
\mathcal{P}(\mathbb{R}^{d})\times \Omega\rightarrow \mathcal{P}(\mathbb{R}%
^{d})$, $(t,\mu ,\omega )\mapsto \lambda(t,\omega )\mu$ such that $\lambda
(0,\omega )=I$, namely the identity map on $\mathcal{P}(\mathbb{R}^{d})$, and 
\begin{equation}
\lambda (t+s,\omega )=\lambda \left( t,\vartheta _{s}\omega \right) \circ
\lambda (s,\omega )  \label{cocycle}
\end{equation}%
for all $t,s\in \lbrack 0,\infty )$ and for all $\omega \in \Omega $. In (%
\ref{cocycle}) the symbol $\circ $ means map composition. $A$ family of maps 
$\lambda (t,\omega )$ satisfying (\ref{cocycle}) is called a cocycle, and (%
\ref{cocycle}) is the cocycle property. The map $\left\{ \vartheta
_{t}:\Omega \rightarrow \Omega \right\} ,t\in \lbrack 0,\infty )$ is a
family of measure-preserving transformations of a probability space $(\Omega
,\mathcal{F},P)$ termed the shift operators;  see, for
instance, Section 2.5 in Karatzas and Shreve \cite{KaratzasShreve} for
further details on the shift operators.

Using this map, the solution of (\ref{ks0}) can be expressed as 
\begin{equation}\label{eq:rds}
\pi _{t}(\omega )=\lambda (t,\omega )\pi ^{0}.
\end{equation}
Moreover for an arbitrary $\mu \in \mathcal{P}(\mathbb{R}^{d})$, the process 
$\pi^\mu =\{\pi^\mu _{t},t\geq 0\}$, defined as 
\begin{equation*}
\pi^{\mu}_t(\omega):=\lambda(t, \omega) \mu
\end{equation*}
is the solution of the SPDE (\ref{ks0}) with initial condition $%
\mu$. Finally the map $\lambda(t, \omega)$ is a continuous map when we endow 
$\mathcal{P}(\mathbb{R}^{d})$, (or, rather, the set of probability measures
with  second moment) with the topology induced by the Wasserstein
metric.

Since the FA process is infinite-dimensional, as explained in Sec.~\ref{ssec:formulation}, 
its RDS characterization is not immediate. RDS theory is well developed 
for finite-dimensional processes \cite{Arnold.1998}. A
subclass of these questions is settled in a fairly satisfactory manner by
the theory of stochastic flows; see, for instance, Ikeda and Watanabe \cite%
{ikedawatanabe} and Kunita \cite{Kunita.1986}. A stochastic flow needs jointly
continuous dependence of the solutions of the stochastic differential equation 
under consideration on time and on the initial state, except for a set of measure zero. 
This often does not hold for infinite dimensions. Some infinite-dimensional systems do generate a
stochastic flow, others do not. For further definitions of possibly
infinite-dimensional RDSs, as well as for other related
results, we refer to Crauel and Flandoli\cite{Crauel.Flandoli.1994} and Flandoli \cite%
{flandoli3}. The RDS characterization of the FA process is discussed in
{\mg Appendix II.E}. 

In this paper, we show that, despite the possible divergence of the prior
distributions, that is, 
\begin{equation*}
\lim_{t\rightarrow \infty }W_{2}(p^{\mu}_{t},p_{t})=\infty ,
\end{equation*}%
the FA process has a stabilizing effect, in the sense
that it keeps the distance $W_{2}(\pi^{\mu}_{t},\pi _{t})$ uniformly bounded in
expectation. Moreover, in the linear case, it makes the distance $%
W_{2}(\pi^{\mu}_{t},\pi _{t})$ vanish asymptotically. The main results of the
paper are Theorems \ref{nonlinear} and \ref{linear} below.

\section{Main Results} \label{sec:wd}

We introduce now the Wasserstein metric on $\mathcal{P}^{2}(\mathbb{R}^{d})$
defined by \eqref{eq:W2} on 
the set of all probability measures on the collection of Borel sets $\mathcal{B}(\mathbb{R}^{d}) $
that have a finite second moment. Recall that
the set ${\Gamma (\mu ,\nu )}$ in \eqref{eq:W2} denotes the collection of all measures on ${{%
\mathbb{R}}^{d}}\times {{\mathbb{R}}^{d}}$ with marginals $\mu $ and $\nu $
on the first and second factor, respectively; it 
is called the set of all couplings of the measures $\mu $ and $\nu $. 

The Wasserstein metric is equivalently defined by 
\begin{equation}  \label{equiv}
{W_{2}(\mu ,\nu )=\left( \inf \mathbb{E}{\big [}|X-Y|^{2}{\big ]}\right)
^{1/2},}
\end{equation}%
where $\mathbb{E}{[Z]}$ denotes the expected value of a random
variable or vector $Z$ and the infimum is taken over all joint distributions of
the random variables $X$ and $Y$ with marginals $\mu $ and $\nu $
respectively.

The main results of the paper are:

\begin{theorem}\label{nonlinear}
For nonlinear coefficients $f$, $h$ and $\sigma$ and measures $\pi_0$ and 
 $\mu\in\mathcal{P}_{2}(\mathbb{R}^{d})$ that satisfy the conditions stated in  {\mg Appendix~II, Sec. B},  
 there exists a bound $R=R\left( \pi_0,\mu\right)$ such that 
\begin{equation} \label{eq:Thm1}
\sup_{t\ge 0}\mathbb E[W_{2}(\pi _{t}^{\mu },\pi_{t})]\leq R.
\end{equation}
\end{theorem}

The complete proof of this theorem is given in {\mg Appendix~II, Sec. B}. For the benefit of the curious but hasty reader, we provide here a brief sketch of the argument. First, we give a bound on the difference $%
\hat{\pi}_{t}^{\mu }-\hat{\pi}_{t}$ between the mean of the FA process initialized from $\mu $ and the original 
FA process, initialized from $\pi _{0}$. This is done in two steps: For arbitrary $\delta>0$, 
we deduce that there exists a constant $c_{\delta }$ independent of $k$ such that 
\begin{equation*}
\sup_{t\in \lbrack k\delta ,\left( k+1\right) \delta ]}\mathbb{E}\left[
\left\vert \hat{\pi}_{t}^{\mu }-\hat{\pi}_{t}\right\vert \right] \leq
c_{\delta }\mathbb{E}\left[ \left\vert \hat{\pi}_{k\delta }^{\mu }-\hat{x}%
_{k\delta }\right\vert \right] . 
\end{equation*}%
This inequality appears as Eq.~\eqref{cd} in the appendix.

Next, we show that there exists $R_{\delta }$ such that $\sup_{k\ge 0}\mathbb{E}\left[ \left\vert \hat{\pi}_{k\delta }^{\mu }-\hat{x}%
_{k\delta }\right\vert \right] \leq R_{\delta }.  $ These two inequalities give us a uniform bound, over $t\in [0,\infty)$, of the difference $\mathbb{E}\left[
\left\vert \hat{\pi}_{t}^{\mu }-\hat{\pi}_{t}\right\vert \right]$. Finally the uniform bound, for all positive times, of $\mathbb{E}[W_{2}(\pi _{t}^{\mu },\pi_{t})] $ comes by means of Lemma \ref{mean} in the appendix
from the bound on the difference $\mathbb{E}\left[
\left\vert \hat{\pi}_{t}^{\mu }-\hat{\pi}_{t}\right\vert \right]$ and that of the  covariance matrices of the
FA process initialized from $\mu $ and, respectively, the original 
FA process, initialized from $\pi _{0}$.

\begin{theorem}\label{linear}
 For linear coefficients $f$, $h$ and $\sigma$ and measures $\pi_0$ and 
 $\mu\in\mathcal{P}_{2}(\mathbb{R}^{d})$ that satisfy the conditions stated in {\mg Appendix~II, Sec. C}, 
 we have the much stronger result that

\begin{equation} \label{eq:Thm2}
\lim_{t\rightarrow \infty }
W_{2}(\pi _{t}^{\mu },\pi_{t})=0.
\end{equation}

\end{theorem}

The complete proof of Theorems \ref{linear} is 
given in {\mg Appendix~II, C}. Again, we provide here a brief sketch of the argument. 
We define $\pi _{t}^{0}:=N\left( \hat{x}_{t}^{0},P_{t}^{0}\right)$
to be a suitably chosen probability measure process that serves as a "reference."  
Then we show that $\pi ^{\mu }$ gets asymptotically close to the
reference process $\pi ^{0}$, regardless of the initial condition, and
since this holds true for $\mu =\pi _{0}$ too, we immediately deduce that $\lim_{t\rightarrow \infty }W_{2}(\pi _{t}^{\mu },\pi _{t})=0$.

The reference process $\pi ^{0}$ is convenient to work with: its centered version  $N\left( 0,P_{t}^{0}\right) $ converges weakly, as well as in Wasserstein distance, to $\pi ^{\infty }=N\left( 0,P_{\infty }\right) $. 
Using the equivalent definition of the Wasserstein distance 
\eqref{equiv}, we deduce that $\lim_{t\rightarrow \infty}W_{2}(\pi _{t}^{\mu },\pi _{t}^{0})=0$, 
if and only if the following three properties hold true:

\begin{itemize}
\item $\lim_{t\rightarrow \infty }| \hat{\pi}_{t}^{\mu }-\hat{x}%
_{t}^{0}| =0;$

\item $\lim_{t\rightarrow \infty }\left\vert P_{\pi _{t}^{\mu
}}-P_{t}^{0}\right\vert =0;$ and

\item $\lim_{t\rightarrow \infty }\left\vert \pi _{t}^{\mu }\left( \varphi
_{t}\right) -\pi _{t}^{0}\left( \varphi _{t}\right) \right\vert =0$ for any
bounded uniformly continuous function $\varphi $, where $\varphi
_{t}$ is the same function shifted by the mean $\hat{x}_{t}^{0},$ that is $%
\varphi _{t}\left( x\right) :=\varphi _{t}\left( x+\hat{x}_{t}^{0}\right) ,$ 
$x\in \mathbb{R}^{d}$. 
\end{itemize}
These three properties are then shown to hold, thus completing the proof. 


\begin{remark} \label{rem:KF}
        In the linear case of Theorem~\ref{linear}, it is fairly easy to verify that the Kalman-Bucy filter \cite{Kalman.1960, Kalman.Bucy.1961} satisfies the assumptions of the theorem. For the nonlinear case of Theorem~\ref{nonlinear}, it is the subject of future research to find numerical criteria that guarantee the required assumptions.   This is a challenging  problem as there are many suboptimal approximations of an optimal filter and the verification of the corresponding hypotheses may prove more difficult. 
\end{remark}

\begin{remark} \label{rem:EKF}
Stability properties of suboptimal filters --- e.g., the extended Kalman-Bucy filters (EKFs) \cite{Ghil.Mal.1991, mt1}, as opposed to the truly optimal filter studied herein  --- have also been studied under the assumption of uniformly stable and fully observable signals. The stability constraint for EKFs has been removed in \cite{bm,bm2,bm3}; see also \cite{km} for a study of the stability of the mean-squared filtering error. 
\end{remark}

\section{Conclusions and Further Work}\label{sec:concl}

\subsection{Summary}\label{ssec:summary}

The main results of this work are given by Theorems~1 and 2 in Sec.~\ref{sec:wd}. Essentially,
\begin{enumerate}[(i)]
        \item 
        For nonlinear dynamics or observations --- including unstable dynamics of the prior process and subject to certain  technical but plausible assumptions --- the supremum of the  expectation of the  Wasserstein distance $W_2$ between the true posterior solution and a solution of the FA process with the wrong initial conditions remains bounded at all future times; and
        \item For linear dynamics and observations --- including unstable dynamics of the prior process and subject to certain technical but plausible assumptions --- the Wasserstein distance $W_2$ between the two posterior distributions tends to zero.
 \end{enumerate}


\subsection{Discussion}\label{ssec:discuss}

These results, to the best of our knowledge, are the first to address the stability of the posterior FA process given an unstable prior process.  In the linear case, convergence in $W_2$ of the posterior processes starting from correct and incorrect initial data, $p_0$ and $\mu$, has been demonstrated (Theorem~\ref{linear}) and the applicability to partial observations has been illustrated in {\mg Appendix~II.D}. In particular, observing the unstable components of the prior process, as originally proposed by A. Trevisan and her collaborators \cite{Carrassi.ea.2008,Carrassi.ea.2007,CTDTU.2008,Trevisan.Uboldi.2004,Trevisan.ea.2010}, seems to be an excellent idea. The results in the nonlinear case only guarantee $W_2$-boundedness of the difference between the  two posterior processes starting from distinct initial probability measures. A considerable amount of practical DA work also indicates that the FA process can track the correct solution \cite{Asch.Bocquet.2016, Bengtsson.ea.1981, Ghil.Mal.1991, Leeuwen.ea.2015}, in particular when using observations from the unstable subspace \cite{Carrassietal2020}.

The first results herein towards a more realistic mathematical treatment of the unstable-dynamics case open the door to a whole slew of additional results, both theoretical and practical.

\paragraph{Deterministic and stochastic EnKF.} Operational DA in NWP relies these days more and more on the EnKF \cite{Bocquet.ea.2010, Carrassietal2020, Leeuwen.2009}. But, in practice, most operational EnKF algorithms  randomize only over the NWP model's initial states and not over observations, too, as done herein.

The operational NWP literature on DA distinguishes, in fact, between the \emph{deterministic} EnKF, which only takes into account random errors in the observations via the covariance matrix $\mathbf R$ of observational errors, and the \emph{stochastic} EnKF, which explicitly simulates random errors in the observations. An excellent review of the EnKF for atmospheric DA in general appears in Houtekamer and Zhang \cite{Houtekamer.2016}, with particular attention to this issue in its Section~2b. Lawson and Hansen \cite{Lawson.Hansen.2004} give interesting examples of the two versions of EnKF being applied to relatively simple examples of atmospheric and oceanic flows in one and two spatial dimensions, and Hoteit et al. \cite{Hoteit.ea.2015} discuss some of the problems that might arise in the {stochastic} EnKF by introducing these observational random errors into the FA process. P.J. van Leeuwen \cite{Leeuwen.2020} has proposed recently a self-consistent way of applying the stochastic EnKF.

Given the novel convergence results obtained herein in the presence of a random observation process and some of the renascent interest in the NWP literature, it might be worthwhile revisiting the usefulness of the stochastic EnKF. In particular, retaining random perturbations in the observations might obviate the need for artificial inflation of the ensemble's rapidly lost variance for the deterministic EnKF.

\paragraph{Multiple models and model error.} In practice, in NWP and elsewhere, prediction can be served by more than one model. The models can differ by their spatial resolution, by the physical processes taken into consideration and by the numerical discretization of the PDEs governing them. Multi-model DA is discussed in some detail by Bach and Ghil \cite{Bach.Ghil.2022}, including the issue of model error growth in this situation. It would be of considerable interest to extend the rigorous results herein to such a broader setting.

Typically, given the fact that DA is more expensive than straight forecasting \cite{Asch.Bocquet.2016, bc, Bengtsson.ea.1981, Carrassietal2020, Ghil.Mal.1991}, it is natural to use lower-resolution models for the FA process than for the forecasting. In the set-up of Sec.~\ref{ssec:setting} herein, doing so corresponds to distinct prior processes $p_t$ and posterior processes $\pi_t$ and appropriate consideration of such issues would be quite worthwhile.


        
\paragraph{Parameter estimation and the synchronization point of view.} In Sec.~\ref{ssec:background}, we have mentioned already the view of the FA process as the synchronization of the forecast model with the observed process \cite{Abarbanel.ea.2017, Duane.ea.2006}. This point of view has been used as a unifying principle between DA and supermodelling, namely the use of ensembles of models that do not only serve for a posteriori averaging of their results but learn from each other in the process of a simulation or prediction run \cite{Duane.ea.2017}. This learning is clearly related to the estimation of imperfectly known model parameters \cite{Gelb.1974,Jazwinski.1970}.
        
In this broader perspective, one could try to demonstrate, given suitable hypotheses and observations, (i) the convergence of a single model's parameter estimation process; and (ii) the convergence of a supermodel to the observed process.

\paragraph{Practical examples.} We presented in {\mg Appendix II.D,}  a simple linear model and two observation schemes to illustrate the fact that stabilization by the FA process does not require observing all of a model's degrees of freedom and that observing just the unstable ones suffices. In future work, we aim to apply these results to realistic models, for example models that numerically approximate the PDEs of geophysical fluid dynamics.   


\paragraph{Particle filters.} The results presented in this paper are theoretical in nature. In practice, the FA process cannot be computed exactly: numerical approximations are required to estimate the posterior distribution of the signal given the data. Among these numerical approximations, particle filters have the crucial property of being theoretically justified in the sense that the numerical error can be controlled by the computational effort. Moreover, they are asymptotically consistent, i.e., as the number of data points used increases, the sequence of estimates converges in probability

A \textit{particle filter} is a sequential Monte Carlo method in which the posterior distributiuon is approximated using a set of \textit{particles}, yielding a random measures of the form
$\displaystyle\sum_{\ell} \mathrm{w}_{t}^{\ell}\delta({x_{t}^{\ell}})$,
where $\delta$ is the Dirac delta function, $\mathrm{w}_{t}^1, \mathrm{w}_{t}^2, \ldots $ are the \textit{weights} of the  particles and $x_{t}^1, x_{t}^2, \ldots$ are their corresponding positions,\cite{bc}
 centered around the state vector $x_{t}$. The approximations evolve in time, by following the time evolution induced by the prior model, and are corrected by the observations to keep them close to the evolution of the FA process \cite{bc, rc}.

Particle filters have been very successful in many applications, including engineering, economics and finance; see, for instance, Doucet et al. \cite{SMCM} and the references therein. In recent years, applications of particle filtering to DA problems for planetary flows have flourished. For in-depth reviews of the most recent efforts in this direction, see refs. \cite{Leeuwen.ea.2015, VetraCarvalho-Leeuwen2018}. Such applications require enhancements of the classical particle filters in order to eable them to tackle the so-called curse of dimensionality, e.g., by relying on optimal transport ideas \cite{Leeuwen.ea.2015}, tempering \cite{CrisanBeskosJasra, KantasBeskosJasra, Wei2, Wei1}, localization \cite{Potthast, Leeuwen.ea.2015}, model reduction \cite{Wei2, Wei3},  jittering\cite{CrisanBeskosJasra}, nudging \cite{Wei1}, and judicious proposal densities \cite{Leeuwen2010}. Some of these approaches have been tested in operational NWP weather prediction systems \cite{Potthast}.  The suitability of particle  filters for high-dimensional problems has been studied in \cite{CrisanBeskosJasra} and tested in \cite{Wei2, Wei1}. For example, in \cite{Wei1}, the method is used for the stochastic incompressible two-dimensional Euler model with forcing and damping, while in \cite{Wei2} it is tested for a two-dimensional quasi-geostrophic model. 

The theoretical results herein pave the way for the stability analysis of particle 
filters under the same assumptions. More precisely, one can attempt to show that particle
filters have numerical errors that can be controlled \emph{uniformly} in time. Again, we will be guided in  pursuing such results by existing ones in nonlinear filtering, stochastic analysis and applied probability \cite{uc1, uc2, uc3, uc4, uc5, uc6}.    The bound will be in expectation, as in Theorem \ref{nonlinear} of this paper. Coupled with Remark \ref{r3} in {\mg Appendix~II, Sec. B}, such results will offer theoretical validation to applying particle filters for \emph{long-run} DA problems.

\newpage

\section*{Appendix I. Forecast Error Growth in NWP}

In order to better understand the nature and role of forecast error growth
in the FA process, we consider here three different models for error growth in
NWP, namely those of C. E. Leith \cite{Leith.1978}, E. N. Lorenz \cite%
{Lorenz.1982} and of A. Dalcher and E. Kalnay \cite{Dalcher.Kalnay.1987}.

The forecast error model of Leith\cite{Leith.1978} is 
\begin{equation}  \label{eq:Leith}
\dot V = \alpha V + S,
\end{equation}
where $V$ is the mean-square error, $S$ is the systematic model error, and $t
$ is the lead time. The $V$ in this appendix should not be confused with the
Brownian motion process in Eq.~\eqref{signal} of Sec.~\ref{sec:framework}.
The forecast error growth in the Leith model is given by 
\begin{equation}
V(t) = \left(V_0 + \frac{S}{\alpha}\right) e^{\alpha t} - \frac{S}{\alpha},
\end{equation}
where $V_0 = V(0)$ is the initial error. Note that short-time forecast
errors grow exponentially and that the systematic model error acts to
increase the coefficient of this growth.

Leith's forecast error model can only apply for short-time error growth,
since it does not saturate. Note that under certain statistical assumptions,
the mean-square error saturation value of a single forecast will be 2$C$,
where $C$ is the climatological variance. In a real NWP model with $N$
variables, the scalar $C$ will be equal to the mean trace of the
climatological covariance matrix $\mathbf{C}$, where $\mathbf{C}$ has
dimension $N \times N$. For ensemble forecasts, the saturation value becomes 
$(1 + 1/m) C$, where $m$ is the ensemble size, see Leith1974\cite{Leith.1974} 

Lorenz's model of forecast error growth \cite{Lorenz.1982} is 
\begin{equation}  \label{eq:Lorenz}
\dot E = aE(E_\infty - E),
\end{equation}
where $E$ is the root-mean-square error and $E_\infty$ is its saturation
value. To compare this model directly to models based on mean-square error,
like Eq.~\eqref{eq:Leith}, we can change variables to $V = E^2$, and get
that 
\begin{equation}  \label{eq:Lorenz_V}
\dot V = 2a V^{1/2}_{\infty} V \left(1 - (V/V_{\infty})^{1/2} \right).
\end{equation}
Lorenz's model includes a nonlinear saturation term, but does not
incorporate systematic model error $S$.

The error model proposed by Dalcher and Kalnay \cite{Dalcher.Kalnay.1987}
(henceforth DK) combines the key features of the Leith \cite{Leith.1974} and
Lorenz \cite{Lorenz.1982} models, 
\begin{equation}  \label{eq:DK}
\dot V = (\alpha V + S)(1 - V/V_\infty);
\end{equation}
it thus includes both saturation $V_\infty$ and systematic model error $S$.
For short-time error growth, we can take $V_\infty\to \infty$, recovering
Leith's model. For $S = 0$, the model is similar to that of Lorenz, but with 
$V_{\infty}$ having unit power in the saturation term, rather than $1/2$.


To compare the three models graphically, we set $\alpha = 1$, $V_0 = 1$, and 
$V_\infty = 100$. To match the Lorenz model's short-term error growth to
that of the other two models, we set $a = \alpha/\left(2
V_{\infty}^{1/2}\right)$.

\begin{figure}[tbp]
\includegraphics[scale=0.5]{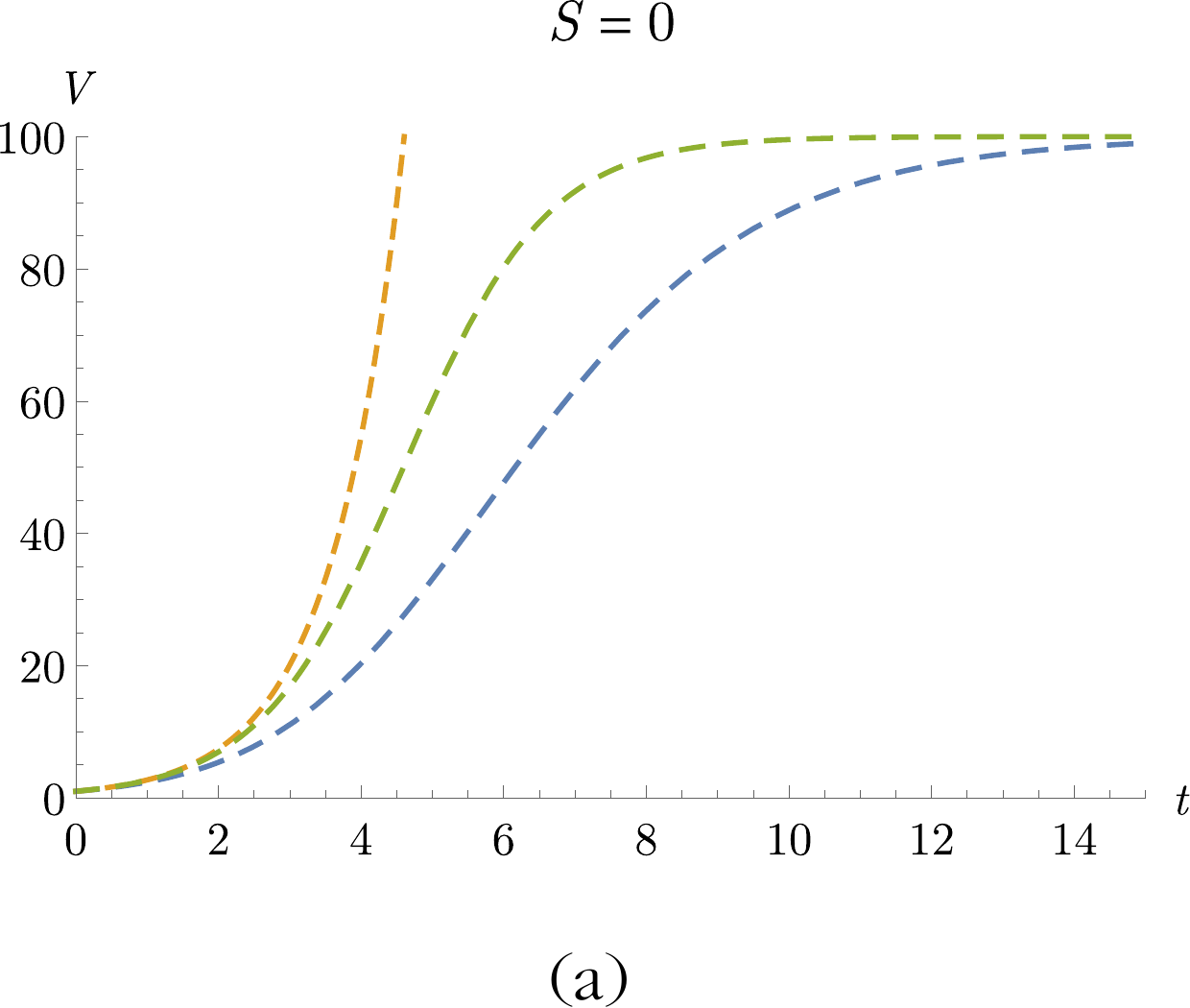} %
\includegraphics[scale=0.5]{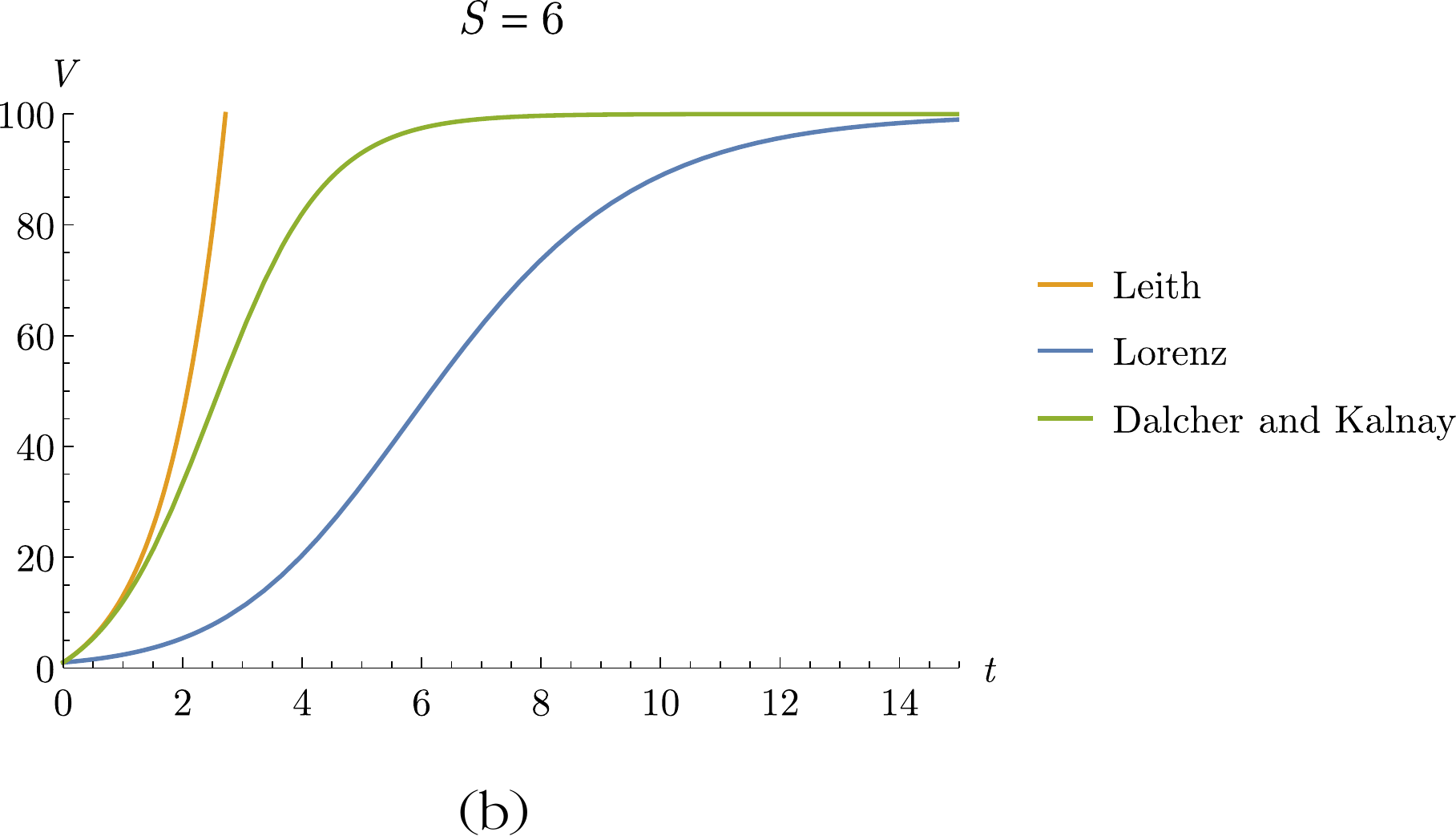}
\caption{Comparison between the model error growth $V(t)$ in the Leith 
\protect\cite{Leith.1978}, Lorenz \protect\cite{Lorenz.1982} and DK 
\protect\cite{Dalcher.Kalnay.1987} models. The systematic model error is (a) 
$S=0$, dashed curves; and (b) $S = 6$, solid curves. See legend for color
identification. }
\label{fig:error_models}
\end{figure}

Figure~\ref{fig:error_models} shows a comparison between the three error
models, both with a perfect model for which $S=0$ and with an imperfect
model with $S=6$. The Leith curve is the same in both cases, since it does
not account for $S$. For the perfect model case, all three curves experience
similar short-term exponential error growth. However, the Leith and DK
models grow faster in the imperfect case than in the perfect case. In both
cases, the Leith model diverges from the other two curves in the medium
range, due to its lack of saturation. The Lorenz and DK curves both saturate
to $V_\infty$, although DK saturates more quickly due to the difference in
the functional forms of the saturation terms in the two models.

DK actually used their model's three parameters, $(\alpha, S, V_\infty)$, to
match their error growth curve to the archived real-time performance for the
years 1980-1981 of the European Centre for Medium-range Weather Forecasts
(ECMWF) model. Their Figure 9 shows separate, near-perfect fits out to 10
days, for boreal winter, when the weather is more active, and boreal summer,
when it is less so. Stroe and Royer \cite{Royer.1993} subsequently
generalized the DK model, by introducing the power $V^p$, with $p \neq 2$,
in the saturation term, and took the limit $p \to \infty$ to obtain 
\begin{equation}  \label{eq:Royer}
\dot V = -aV \log(V/V_\infty).
\end{equation}
These authors found that Eq.~\eqref{eq:Royer} gave better fits for
extended-range, 45-day experimental weather forecasts than either Eq.~%
\eqref{eq:Lorenz_V} or \eqref{eq:DK}.

Simmons et al. \cite{Simmons.ea.1995} also obtained rather good fits to
operational NWP model performance with Lorenz's quadratic error growth
model, but with a smaller error growth exponent that in the DK paper \cite%
{Dalcher.Kalnay.1987}. Trevisan et al. \cite{Trevisan.ea.1992}, though,
showed --- by using an intermediate, quasigeostrophic two-layer model on the 
$\beta$-plane \cite{Malguzzi.ea.1990} --- that only very small initial
errors in such a model obey Lorenz's quadratic error growth model and that
the error growth curve in general depends significantly on the magnitude of
the initial errors.

Savij\"arvi \cite{Savijarvi.1995} combined features of the Lorenz and DK
models in the study of the (then) U.S. National Meteorological Center's
(NMC's) Medium-Range Forecast (MRF) Model's 0--10-day forecasts for
1988--1993. Growth parameters, as well as model and analysis errors for this
data set, were estimated using the quadratic error growth assumption.
Savij\"arvi showed that both the MRF model error and analysis error nearly
halved during the six years under study but, at the same time, the growth
parameters nearly doubled, since smaller errors grow faster.

Model error growth is thus a complex topic with much more to be said about;
see, for instance, the line of inquiry developed by C. Nicolis and coworkers, which
includes transient bimodality of the error's pdf \cite{Nicolis.ea.2009}. %
The topic's quick review in this appendix suffices, though, to 
show the presence of error growth--generating instabilities in high-end,
operational NWP models. The ground covered here in Secs.~\ref*%
{sec:framework} and \ref*{sec:wd} and in Appendix~II shows that DA can 
overcome these instabilities,  in theory as well as in practice.

\section*{Appendix II. Rigorous Definitions and Proofs}

\subsection*{A. The Wasserstein topology}

In this appendix, we present the reason for showing the stabilizing effect of the FA process on unstable dynamics with respect to the Wasserstein topology --- i.e., the topology generated by
the Wasserstein distance --- and not with respect to the more popular weak
topology. Before doing so, we provide a little more information on the 
Wasserstein distance for the benefit of NWP practitioners who might not yet
be familiar with it.

Gaspard Monge, an artillery officer in Napoleon's armies, as well as one of the founders of France's Ecole Polytechnique, introduced it as early as the 1780s \cite{Monge.1781}, and Leonid V. Kantorovich\cite{Kantorovich.1942} used it during World War II in optimizing the transport of resources within the Soviet Union. The contemporary developments of this distance and of its applications are largely due to R.~L.~Dobrushin \cite{Dobrushin.1970}, who coined the name  Wasserstein \cite{Wasserstein.1969} distance for it, and to C. Villani \cite{Villani.2009}. 

In the climate sciences, Ghil \cite{Ghil.2015} illustrated the use of the Wasserstein distance for measuring the parameter sensitivity of simple models with time-dependent forcing, thus providing a link between nonautonomous dynamical systems theory \cite{Caraballo.Han.2017} and optimal transport \cite{Villani.2009}. Robin et al. \cite{Robin.ea.2017} then used this distance to compute the difference between the snapshot attractors of the Lorenz \cite{Lorenz.1984} model for different time-dependent forcings, while Vissio et al. \cite{Vissio.ea.2020} used it to help intercomparing climate models and evaluating their performance against given benchmarks in the Coupled Model Intercomparison Project that is part of the  Intergovernmental Panel on Climate Change process.

Returning to the main purpose of this appendix, let us introduce first some notation. For a measure $%
\mu\in \mathcal{P}^2({\mathbb{R}^d})$, we use the following notation:
\begin{enumerate}[(i)]
\item {\bf means}: $\hat{\mu}=(\hat{\mu}^{i})_{i\in 1,\dots d}$ is the mean vector of $\mu$, i.e.,
\begin{equation} \label{eq:means}
\hat{\mu^{i}}=\int_{\mathbb{R}^d}x_{i}\mu(\mathrm{d }x), \quad {i\in
1,\dots d}; \quad|\hat{\mu}| = \left(\sum_{i=1}^d(\hat{\mu^{i}})^2\right)^{1/2};
\end{equation}

\item {\bf second moments}: $\mu^2$ is the sum of the second moments of $\mu$,  i.e.,
\begin{equation} \label{eq:moments}
        \mu^{2}=\sum_{i=1}^d\int_{\mathbb{R}^d}x_{i}^2\mu(\mathrm{d }x); 
\end{equation}

\item {\bf covariance matrix}: $P_\mu=(P_\mu^{ij})_{i,j\in 1,\dots d}$ is the covariance matrix of $%
\mu$, i.e., 
\begin{subequations} \label{eq:covar}
        \begin{align}
& P_\mu^{ij} = \int_{\mathbb{R}^d}(x_i-\hat{\mu^{i}})(x_j-\hat{\mu^{j}})\mu(%
\mathrm{d }x), \ \ {i,j\in 1,\dots d,} \\
& |P_\mu| = \left(\sum_{i,j=1}^d(P_\mu^{ij})^2 \right)^{1/2}.
        \end{align}
\end{subequations}
\end{enumerate}
Next, we recall the definition of the weak topology on the space of probability
measures $\mathcal{P}(\mathbb{R}^{d})$ and note that,
of course, $\mathcal{P}^2(\mathbb{R}^{d})\subset \mathcal{P}(\mathbb{R}^{d})$.
We can thus consider also weak convergence of probability measures belonging to
the smaller space $\mathcal{P}^2(\mathbb{R}^{d})$:

\begin{definition}[Weak topology]
\label{defn:condExpect:weakConv} A sequence of probability measures $\left(
\mu _{n}\right) _{n}\in \mathcal{P}(\mathbb{R}^{d})$, converges \emph{weakly}
to $\mu \in \mathcal{P}(\mathbb{R}^{d})$ if and only if $\left( \mu
_{n}\left( \varphi \right) \right) _{n}$ converges to $\mu \left( \varphi
\right) $ as $n\rightarrow \infty $ for all $\varphi \in C_{b}(\mathbb{R}%
^{d})$. The {\em weak topology} on the space $\mathcal{P}(\mathbb{R}^{d})$ is
defined to be the weakest topology such that for all $f\in C_{b}(\mathbb{R}%
^{d})$, the function $\mu \mapsto \mu \left( f\right) $ is continuous.
\end{definition}

The weak convergence of $\left( \mu _{n}\right) _{n}$ to $\mu $ is denoted $%
\mu _{n}\Rightarrow \mu $. A set of probability measures $\mathcal{A}\subset 
\mathcal{P}(\mathbb{R}^{d})$ is relatively compact in the weak topology if and only if
for all $\varepsilon >0$ there exists $K_{\varepsilon }$ such that $\mu
\left( K_{\varepsilon }\right) \geq 1-\varepsilon $\thinspace\ for all $\mu
\in \mathcal{A}$. If $\mathcal{A}\subset \mathcal{P}^2(\mathbb{R}^{d})$, one
can show that the set $\mathcal{A}$ will be relatively compact in the weak
topology if the means and the covariance matrices of the probability
measures in the set $\mathcal{A}$ are uniformly bounded. 

The set $\mathcal{A}$ can, however, still be relatively compact, even if the means of the
probability measures in it do not remain bounded. For
example, if we choose $\mu _{n}=\left( 1-\frac{1}{n}\right) \delta _{0}+%
\frac{1}{n}\delta _{n^{2}}, $ then the sequence $\mu _{n}$ is relatively
compact --- in fact, $\mu _{n}\Rightarrow \delta_0 $ --- but the
corresponding sequence of means $\hat \mu _{n} $ is not bounded, since $\mu _{n}\left( \varphi
\right) =n.$ On the other hand, if we choose $\mu _{n}=\left( 1-\frac{1}{n}\right)
\delta _{0}+\frac{1}{n}\delta _{n}, $ then the sequence $\mu _{n}$ is
relatively compact, the means $\hat\mu\equiv1$ form a
trivially bounded sequence, $\mu _{n}\left( \varphi \right) =1$, but the
second moments are not, as $\mu _{n}^2=n$. This state of affairs is not
satisfactory for our purposes. 

The Wasserstein topology, though, adds the convergence of the first and second moments to
the weak convergence of the measures. To be precise, we have $%
\lim_{n\rightarrow \infty }{W_{2}({\mu }_{n},\mu )=0}$ for $\mu_{n},\mu \in 
\mathcal{P}^2({\mathbb{R}^d})$, if and only if ${\mu }_{n} $ converges to $%
\mu $ in the weak topology and the first and second moments converge as well,
using the notation of Eqs.~\eqref{eq:means}--\eqref{eq:covar} above. 
Moreover, a set of probability measures $\mathcal{%
A}\in \mathcal{P}^2({\mathbb{R}^d})$ is relatively compact in the topology
given by the Wasserstein distance if and only if 
\begin{equation}
\lim_{R\mapsto\infty}\sup_{\mu\in\mathcal{A}}\int_{|x|>R}|x|^2\mu(\mathrm{d }%
x)=0.
\end{equation}
Finally, we have the following lemma which follows immediately from the alternative
definition \eqref{equiv} of the Wasserstein distance $W_2$: 
\begin{lemma}\label{mean} There exists a constant $C=C(d)$ such that, for any  $\mu ,\nu\in \mathcal{P}^2({\mathbb{R}^d})$, 
\begin{eqnarray*}
(W_{2}(\mu ,\nu ))^2&\le& C (\mu^2+ \nu^2 ),\\ 
W_{2}(\mu ,\nu )&\le& C(|P_\mu|^{\frac{1}{2}}+|P_\nu|^{\frac{1}{2}}
+|\hat\mu-\hat\nu|).
\end{eqnarray*}

\end{lemma}

\subsection*{B. Assumptions and proof of Theorem~\ref{nonlinear}}

\label{aps}

To start, we formulate here the set of assumptions on the coefficients of
the signal and observation equations \eqref{signal}+\eqref{observation}
under which Theorem \ref{nonlinear} holds:

\begin{itemize}
\item We assume that the coefficients $f$ and $h$ can be decomposed into a
linear part and a \emph{bounded} nonlinear part. 
In other words, we will assume that 
\begin{equation}
f=F{\mathcal{I}} +\tilde{f},\ \ h=H{\mathcal{I}} +\tilde{h},  \label{Fhsigma}
\end{equation}%
where 
\begin{enumerate}[(i)]
        \item ${\mathcal{I}} :\mathbb{R}^{d}\rightarrow \mathbb{R}^{d}$ is the
identity function defined as ${\mathcal{I}} (x)=x$ for any $x\in \mathbb{R}%
^{d}$;
\item $F\in \mathbb{R}^{d\times d}$, $H\in \mathbb{R}^{d\times n}$ are
given matrices; and
\item $\tilde{f}:\mathbb{R}^{d}\rightarrow \mathbb{R}^{d}$ , $\tilde{h}:%
\mathbb{R}^{d}\rightarrow \mathbb{R}^{n}$ are bounded measurable functions
that incorporate the nonlinear parts of the coefficients of the system %
\eqref{signal}+\eqref{observation}.
\end{enumerate}

\item The covariance matrices of the processes $\pi ^{\mu }$ and $%
\pi $, respectively, are uniformly bounded in expectation: 
\begin{equation}
\sup_{t\geq 0}\mathbb{E}[|P_{t}^{\mu }|^{8}]=C^{\mu }<\infty ,\ \
\sup_{t\geq 0}\mathbb{E}[|P_{t}^{\pi _{0}}|^{8}]=C^{\pi _{0}}<\infty .
\label{controlpt}
\end{equation}

\item The matrix-valued process $Q_{s}^{\mu }:=F-P_{s}^{\mu }H^{\top
}H-\lambda _{s}^{\mu }H$ is exponentially stable in expectation, where $%
\lambda ^{\mu }$ is the matrix-valued process defined as $\lambda _{s}^{\mu
} := \pi _{s}^{\mu }({\mathcal{I}} -\hat{\pi}^{\mu })\tilde{h}^{\top }$ for 
$s\geq 0$. In other words, if $\psi _{s:t}$ is the solution of the linear
matrix ordinary differential equation 
\begin{equation*}
\mathrm{d}_{t}\psi _{s:t} =Q_{t}^{\mu }\psi _{s:t},\ \ \psi _{s:s}=I, 
\end{equation*}
where $I$ is the identity matrix, then there exist some constant $c>0$ such
that 
\begin{equation}
\mathbb{E}\left[ \left\vert \psi _{s:t}\right\vert ^{2}\right] \leqq
e^{-c\left( t-s\right) }.  \label{exp}
\end{equation}%
Moreover, we assume that there exists $\delta >0$ such that 
\begin{equation}  \label{psi-}
\int_{k\delta }^{\left( k+1\right) \delta }\mathbb{E}\left[ \left\vert \psi
_{k\delta :s}^{-1}\right\vert ^{4}\right] \mathrm{d}s\leq C^{inv}
\end{equation}%
where $C^{inv}$ is a constant independent of $k$.\bigskip \newline
\end{itemize}

The proof of Theorem \ref{nonlinear} requires the bound of the difference $%
\hat{\pi}_{t}^{\mu }-\hat{\pi}_{t}$ between the mean of the
FA process initialized from $\mu $ and the original 
FA process, initialized from $\pi _{0}$. We deduce from 
\eqref{ks0} 
 that, for $t\in \lbrack k\delta ,\left( k+1\right)
\delta ]$, 
\begin{eqnarray*}
\hat{\pi}_{tk\delta }^{\mu } &=&\hat{\pi}_{k\delta }^{\mu }+\int_{k\delta
}^{t}\pi _{s}^{\mu }(f)\,\mathrm{d}s+\int_{k\delta }^{t}\varkappa _{s}^{\mu
}(\mathrm{d}Y_{s}-\pi _{s}^{\mu }(h)\,\mathrm{d}s), \\
\hat{\pi}_{t}^{\mu } &=&\hat{\pi}_{k\delta }^{\mu }+\int_{k\delta }^{t}\pi
_{s}^{\mu }(f)\,\mathrm{d}s+\int_{k\delta }^{t}\varkappa _{s}^{\mu }\left( 
\mathrm{d}I_{s}+\left( \pi _{s}(h)-\pi _{s}^{\mu }(h)\right) \mathrm{d}%
s\right), \\
\hat{\pi}_{t} &=&\hat{\pi}_{k\delta }+\int_{k\delta }^{t}\pi _{s}(f)\,%
\mathrm{d}s+\int_{k\delta }^{t}\varkappa _{s}\mathrm{d}I_{s},
\end{eqnarray*}%
Here%
\begin{eqnarray*}
\varkappa _{s}^{\mu } &:&=\pi _{s}^{\mu }({\mathcal{I}} h^{\top })-\pi
_{s}^{\mu }({\mathcal{I}} )\pi _{s}^{\mu }(h^{\top })=P_{s}^{\mu }H^{\top
}+\lambda _{s}^{\mu }, \\
\varkappa _{s} &:&=\pi _{s}({\mathcal{I}} h^{\top })-\pi _{s}({\mathcal{I}}
)\pi _{s}(h^{\top })=P_{s}^{\pi _{0}}H^{\top }+\lambda _{s},
\end{eqnarray*}%
$P_{s}^{\mu }$, $P_{s}^{\pi _{0}}$ are the covariance matrices of $\pi
_s^{\mu }$ and of $\pi $, respectively, and $\lambda _{s}^{\mu }:=\pi _{s}^{\mu
}(({\mathcal{I}} -\hat{\pi}^{\mu })\tilde{h}^{\top })$ and $%
\lambda _{s}:=\pi _{s}(({\mathcal{I}} -\hat{\pi})\tilde{h}^{\top })$,  respectively. It
follows that 
\begin{eqnarray}
\hat{\pi}_{t}^{\mu }-\hat{\pi}_{t} &=&\left( \hat{\pi}_{k\delta }^{\mu }-%
\hat{\pi}_{k\delta }\right) +\int_{k\delta }^{t}\left( \pi _{s}^{\mu }-\pi
_{s}\right) (f-\varkappa _{s}^{\mu }h)\,\mathrm{d}s+\int_{k\delta
}^{t}\left( \varkappa _{s}^{\mu }-\varkappa _{s}\right) \mathrm{d}I_{s} 
\notag \\
&=&\left( \hat{\pi}_{k\delta }^{\mu }-\hat{\pi}_{k\delta }\right)
+\int_{k\delta }^{t}Q_{s}^{\mu }\left( \hat{\pi}_{s}^{\mu }-\hat{\pi}%
_{s}\right) \mathrm{d}s+\int_{k\delta }^{t}\left( P_{s}^{\mu }-P_{s}\right)
H^{\top }\mathrm{d}I_{s}+z_{k\delta :t}^{\mu },  \label{deltapi}
\end{eqnarray}%
where $z_{k\delta :t}^{\mu }$ is a process that contains the nonlinearities
in the evolutions $\hat{\pi}^{\mu }$ and $\hat{\pi}$%
\begin{equation*}
z_{k\delta :t}^{\mu }=\int_{k\delta }^{t}\left( \pi _{s}^{\mu }-\pi
_{s}\right) (\tilde{f}-\varkappa _{s}^{\mu }\tilde{h})\mathrm{d}%
s+\int_{k\delta }^{t}\left( \lambda _{s}^{\mu }-\lambda _{s}\right) \mathrm{d%
}I_{s},~~~t\in \lbrack k\delta ,\left( k+1\right) \delta ]. 
\end{equation*}%

Replacing the observation process $Y_{t}$ by the innovation process $I_{t}$
in the evolution equations for $\hat{\pi}_{t}^{\mu }$ and $\hat{\pi}%
_{t}$, respectively, is important: Unlike $Y_t$, $I_{t}$ is a standard Brownian motion,
which enables us to use classical stochastic calculus properties to bound
the moments of the stochastic integrals appearing in Eq.~\eqref{deltapi}. 

Using an argument based on the Gr\"onwall inequality, one deduces that there exists a
constant $c_{\delta }$ independent of $k$ such that 
\begin{equation}
\sup_{t\in \lbrack k\delta ,\left( k+1\right) \delta ]}\mathbb{E}\left[
\left\vert \hat{\pi}_{t}^{\mu }-\hat{\pi}_{t}\right\vert \right] \leq
c_{\delta }\mathbb{E}\left[ \left\vert \hat{\pi}_{k\delta }^{\mu }-\hat{x}%
_{k\delta }\right\vert \right] . \label{cd}
\end{equation}%
We will show that there exists $R_{\delta }$ such that 
\begin{equation}
\sup_{k\ge 0}\mathbb{E}\left[ \left\vert \hat{\pi}_{k\delta }^{\mu }-\hat{x}%
_{k\delta }\right\vert \right] \leq R_{\delta }.  \label{rd}
\end{equation}%
From (\ref{cd}) and (\ref{rd}), one can then deduce that 
\begin{eqnarray*}
\sup_{t\geq 0}\mathbb{E}\left[ \left\vert \hat{\pi}_{t}^{\mu }-\hat{\pi}%
_{t}\right\vert \right] &\leq& \sup_{k \ge 0}\sup_{t\in \lbrack k\delta
,\left( k+1\right) \delta ]}\mathbb{E}\left[ \left\vert \hat{\pi}_{t}^{\mu }-%
\hat{\pi}_{t}\right\vert \right] \\
&\leq& c_{\delta }\sup_{k \ge 0}\mathbb{E}%
\left[ \left\vert \hat{\pi}_{k\delta }^{\mu }-\hat{x}_{k\delta }\right\vert %
\right] \leq c_{\delta }R_{\delta }<\infty . 
\end{eqnarray*}%
Finally, from Lemma \ref{mean} and \eqref{controlpt}, it follows that 
\begin{eqnarray*}
\sup_{t\geq 0}\mathbb{E}[W_{2}(\pi _{t}^{\mu },\pi _{t})] &\leq &\sup_{t\geq
0}\mathbb{E}\left[ \sqrt{|P_{t}^{\mu }|}\right] +\sup_{t\geq 0}\mathbb{E}%
\left [\sqrt{|P_{t}^{\pi _{0}}|}\right]+\sup_{t\geq 0}\mathbb{E}\left[
\left\vert \hat{\pi}_{t}^{\mu }-\hat{\pi}_{t}\right\vert \right] \\
&\leq &\left( \sup_{t\geq 0}\mathbb{E}\left[ |P_{t}^{\mu }|^{8}\right]
\right) ^{\frac{1}{16}}+\left( \sup_{t\geq 0}\mathbb{E}[|P_{t}^{\pi
_{0}}|^{8}]\right) ^{\frac{1}{16}}+\sup_{t\geq 0}\mathbb{E}\left[ \left\vert 
\hat{\pi}_{t}^{\mu }-\hat{\pi}_{t}\right\vert \right] \\
&\leq &\left( C^{\mu }\right) ^{\frac{1}{16}}+\left( C^{\pi _{0}}\right) ^{%
\frac{1}{16}}+c_{\delta }R_{\delta }<\infty,
\end{eqnarray*}%
which yields our claim.

To complete the proof it remains to show the validity of \eqref{rd}. From \eqref{deltapi},
one obtains a form of Duhamel's principle for the difference between the values
of the mean $\hat{\pi}$ of the FA process at the steps $\left( k+1\right) \delta$ and
$k \delta$:
\begin{eqnarray}
\hat{\pi}_{\left( k+1\right) \delta }^{\mu }-\hat{\pi}_{k\delta }^{\mu }
&=&\psi _{k\delta ,\left( k+1\right) \delta }\left( \hat{\pi}_{k\delta
}^{\mu }-\hat{x}_{k\delta }\right) +\psi _{k\delta ,\left( k+1\right) \delta
}\int_{k\delta }^{\left( k+1\right) \delta }\psi _{k\delta
,,s}^{-1}(P_{s}^{\mu }-P_{s})H^{\top }\mathrm{d}I_{s}  \notag \\
&&+\psi _{k\delta ,\left( k+1\right) \delta }\int_{k\delta }^{\left(
k+1\right) \delta }\psi _{k\delta ,s}^{-1}\mathrm{d}z_{s}^{\mu }.
\label{deltapi2}
\end{eqnarray}%
We analyze next each of the three terms on the right-hand side of (\ref%
{deltapi2}). For the first term, we use (\ref{exp}) to derive the inequality%
\begin{equation}
\mathbb{E}\left[ \left\vert \psi _{k\delta ,\left( k+1\right) \delta }\left( 
\hat{\pi}_{k\delta }^{\mu }-\hat{x}_{k\delta }\right) \right\vert \right]
\leq e^{-c\delta }\mathbb{E}\left[ \left\vert \hat{\pi}_{k\delta }^{\mu }-%
\hat{x}_{k\delta }\right\vert \right].  \label{ld1}
\end{equation}%
For the second term, we use the so-called It\^{o}'s integral isometry
property, cf.~Karatzas and Shreve\cite{KaratzasShreve}, to obtain
\begin{eqnarray}
&&\hspace{-2.5cm}\mathbb{E}\left[ \left\vert \psi _{k\delta ,\left(
k+1\right) \delta }\int_{k\delta }^{\left( k+1\right) \delta }\psi _{k\delta
,,s}^{-1}(P_{s}^{\mu }-P_{s})H^{\top }\mathrm{d}I_{s}\right\vert \right] ^{2}
\notag \\
&\leq &\mathbb{E}\left[ \left\vert \psi _{k\delta ,\left( k+1\right) \delta
}\right\vert ^{2}\right] \mathbb{E}\left[ \left\vert \int_{k\delta }^{\left(
k+1\right) \delta }\psi _{k\delta ,,s}^{-1}(P_{s}^{\mu }-P_{s})H^{\top }%
\mathrm{d}I_{s}\right\vert ^{2}\right]  \notag \\
&\leq &e^{-2c\delta }\mathbb{E}\left[ \int_{k\delta }^{\left( k+1\right)
\delta }\left\vert \psi _{k\delta ,,s}^{-1}(P_{s}^{\mu }-P_{s})H^{\top
}\right\vert ^{2}\mathrm{d}s\right]  \notag \\
&\leq &e^{-2c\delta }\mathbb{E}\left[ \int_{k\delta }^{\left( k+1\right)
\delta }\left\vert \psi _{k\delta ,,s}^{-1}\right\vert ^{4}\mathrm{d}s\right]
+e^{-2c\delta }\mathbb{E}\left[ \int_{k\delta }^{\left( k+1\right) \delta
}\left\vert (P_{s}^{\mu }-P_{s})H^{\top }\right\vert ^{4}\mathrm{d}s\right] 
\notag \\
&\leq &C^{I}e^{-2c\delta },  \label{ld2}
\end{eqnarray}
where $C^{I}=C^{inv}+\delta \left\vert H\right\vert ^{4}\sqrt{C^{\mu }}%
+\delta \left\vert H\right\vert ^{4}\sqrt{C^{\pi _{0}}}$.

Next, \eqref{controlpt} implies that there exists positive constants $%
C^{\varkappa ^{\mu }},C^{\varkappa },C^{\lambda ^{\mu }},C^{\lambda }$ such
that 
\begin{equation}
\sup_{t\geq 0}\mathbb{E}[|\varkappa _{t}^{\mu }|^{8}]=C^{\varkappa ^{\mu
}}<\infty ,\ \ \sup_{t\geq 0}\mathbb{E}[|\varkappa _{t}|^{8}]=C^{\varkappa
}<\infty ,~~~\sup_{t\geq 0}\mathbb{E}[|\lambda _{t}^{\mu }|^{8}]=C^{\lambda
^{\mu }}<\infty ,~~~\sup_{t\geq 0}\mathbb{E}[|\lambda _{s}|^{8}]=C^{\lambda
}<\infty .  \label{othercontrol}
\end{equation}%
For the third term, one gets%
\begin{eqnarray}
\int_{k\delta }^{\left( k+1\right) \delta }\psi _{k\delta ,s}^{-1}\mathrm{d}%
z_{s}^{\mu }&=&\int_{k\delta }^{\left( k+1\right) \delta }\psi _{k\delta
,s}^{-1}\left( \pi _{s}^{\mu }-\pi _{s}\right) (\tilde{f}-\varkappa
_{s}^{\mu }\tilde{h})\mathrm{d}s  \notag \\
&&+\int_{k\delta }^{\left( k+1\right) \delta }\psi _{k\delta ,s}^{-1}\left(
\pi _{s}^{\mu }-\pi _{s}\right) (\tilde{f}-\varkappa _{s}^{\mu }\tilde{h}%
)\left( \lambda _{s}^{\mu }-\lambda _{s}\right) \mathrm{d}I_{s},  \label{ld3}
\end{eqnarray}%
and one proceeds to bound separately the two terms in (\ref{ld3}). For the first term in (%
\ref{ld3}), we get 
\begin{eqnarray}
&&\hspace{-2cm}\mathbb{E}\left[ \left\vert \psi _{k\delta ,\left( k+1\right)
\delta }\int_{k\delta }^{\left( k+1\right) \delta }\psi _{k\delta
,s}^{-1}\left( \pi _{s}^{\mu }-\pi _{s}\right) (\tilde{f}-\varkappa
_{s}^{\mu }\tilde{h})\mathrm{d}s\right\vert \right]  \notag \\
&\leq &e^{-c\delta }\mathbb{E}\left[ \int_{k\delta }^{\left( k+1\right)
\delta }\left\vert \psi _{k\delta ,s}^{-1}\left( \pi _{s}^{\mu }-\pi
_{s}\right) (\tilde{f}-\varkappa _{s}^{\mu }\tilde{h})\right\vert \mathrm{d}s%
\right]  \notag \\
&\leq &e^{-c\delta }\left( \mathbb{E}\left[ \int_{k\delta }^{\left(
k+1\right) \delta }\left\vert \psi _{k\delta ,s}^{-1}\right\vert ^{2}\mathrm{%
d}s\right] +\mathbb{E}\left[ \int_{k\delta }^{\left( k+1\right) \delta
}\left\vert \left( \pi _{s}^{\mu }-\pi _{s}\right) (\tilde{f}-\varkappa
_{s}^{\mu }\tilde{h})\right\vert ^{2}\mathrm{d}s\right] \right)  \notag \\
&\leq &e^{-c\delta }\left( C^{inv}+\delta +\left\vert \tilde{f}\right\vert
^{2}\delta +\left\vert \tilde{h}\right\vert ^{2}\left( C^{\varkappa ^{\mu
}}\right) ^{\frac{1}{4}}\delta \right) .  \label{ld4}
\end{eqnarray}%
For the second term in (\ref{ld3}), we get that 
\begin{eqnarray}
&&\hspace{-1cm}\mathbb{E}\left[ \left\vert \psi _{k\delta ,\left( k+1\right)
\delta }\int_{k\delta }^{\left( k+1\right) \delta }\psi _{k\delta
,s}^{-1}\left( \pi _{s}^{\mu }-\pi _{s}\right) (\tilde{f}-\varkappa
_{s}^{\mu }\tilde{h})\left( \lambda _{s}^{\mu }-\lambda _{s}\right) \mathrm{d%
}I_{s}\right\vert \right]  \notag \\
&\leq &e^{-c\delta }\mathbb{E}\left[ \left\vert \int_{k\delta }^{\left(
k+1\right) \delta }\psi _{k\delta ,s}^{-1}\left( \pi _{s}^{\mu }-\pi
_{s}\right) (\tilde{f}-\varkappa _{s}^{\mu }\tilde{h})\left( \lambda
_{s}^{\mu }-\lambda _{s}\right) \mathrm{d}I_{s}\right\vert ^{2}\right] 
\notag \\
&\leq &e^{-c\delta }\int_{k\delta }^{\left( k+1\right) \delta }\mathbb{E}%
\left[ \left\vert \psi _{k\delta ,s}^{-1}\left( \pi _{s}^{\mu }-\pi
_{s}\right) (\tilde{f}-\varkappa _{s}^{\mu }\tilde{h})\left( \lambda
_{s}^{\mu }-\lambda _{s}\right) \right\vert ^{2}]\mathrm{d}s\right]  \notag
\\
&\leq &e^{-c\delta }\int_{k\delta }^{\left( k+1\right) \delta }\mathbb{E}%
\left[ \left\vert \psi _{k\delta ,s}^{-1}\right\vert ^{4}]\mathrm{d}s\right]
+e^{-c\delta }\int_{k\delta }^{\left( k+1\right) \delta }\mathbb{E}\left[
\left\vert \left( \pi _{s}^{\mu }-\pi _{s}\right) (\tilde{f}-\varkappa
_{s}^{\mu }\tilde{h})\left( \lambda _{s}^{\mu }-\lambda _{s}\right)
\right\vert ^{4}]\mathrm{d}s\right]  \notag \\
&\leq &e^{-c\delta }\left( C^{inv}+\delta \left\vert \tilde{f}\right\vert
^{4}\left( \sqrt{C^{\lambda ^{\mu }}}+\sqrt{C^{\lambda }}\right) +\delta
\left\vert \tilde{h}\right\vert ^{4}\left( C^{\varkappa ^{\mu
}}+C^{\varkappa }\right) \right) .  \label{ld5}
\end{eqnarray}%
From (\ref{ld3}), (\ref{ld4}) and (\ref{ld5}), it follows that there exists a
constant 
\begin{equation*}
C^{z^{\mu }}=C^{z^{\mu }}\left( \delta ,C^{inv},\tilde{f},C^{\lambda ^{\mu
}},C^{\lambda },\tilde{h},C^{\varkappa ^{\mu }},C^{\varkappa }\right) 
\end{equation*}
independent of $k$ such that 
\begin{equation}
\mathbb{E}\left[ \left\vert \int_{k\delta }^{\left( k+1\right) \delta }\psi
_{k\delta ,s}^{-1}\mathrm{d}z_{s}^{\mu }\right\vert \right] \leq C^{z^{\mu
}}e^{-c\delta } .  \label{ld6}
\end{equation}%
Finally, from (\ref{deltapi2}), (\ref{ld1}), (\ref{ld2}) and (\ref{ld6}) we
deduce that 
\begin{equation}
\mathbb{E}\left[ \left\vert \hat{\pi}_{\left( k+1\right) \delta }^{\mu }-%
\hat{x}_{\left( k+1\right) \delta }\right\vert \right] \leq e^{-c\delta
}\left( \mathbb{E}\left[ \left\vert \hat{\pi}_{k\delta }^{\mu }-\hat{x}%
_{k\delta }\right\vert \right] +C^{I}+C^{z^{\mu }}\right) .  \label{ld7}
\end{equation}%

Choose now $R=R\left( \delta \right) =\max \left( \left\vert \hat{\mu}-\hat{%
\pi}_{0}\right\vert ,\frac{C^{I}+C^{z^{\mu }}}{c\delta }\right) $ and use
induction to prove (\ref{rd}). From the definition of $R$, we deduce that $%
\left\vert \hat{\mu}-\hat{\pi}_{0}\right\vert \leq R$. Next assume that $%
\mathbb{E}\left[ \left\vert \hat{\pi}_{k\delta }^{\mu }-\hat{x}_{k\delta
}\right\vert \right] \leq R$. From this and (\ref{ld7}), one can obtain that 
\begin{equation*}
\mathbb{E}\left[ \left\vert \hat{\pi}_{\left( k+1\right) \delta }^{\mu }-%
\hat{x}_{\left( k+1\right) \delta }\right\vert \right] \leq e^{-c\delta
}\left( R+C^{I}+C^{z^{\mu }}\right) \leq e^{-c\delta }R\left( 1+\frac{%
C^{I}+C^{z^{\mu }}}{R}\right) \leq e^{-c\delta }R\left( 1+c\delta \right)
\leq R. 
\end{equation*}%
It follows that $\mathbb{E}\left[ \left\vert \hat{\pi}_{k\delta }^{\mu }-%
\hat{x}_{k\delta }\right\vert \right] \leq R$ holds true for any $k \ge 0$
and so does (\ref{rd}). The proof of Theorem \ref{nonlinear} is now
complete. \qed

\begin{remark}\label{power8} 

The proof of the Theorem \ref{nonlinear} relies on the application of Duhamel's principle, as in (\ref{deltapi2}). The Duhamel principle can be applied on an arbitrary interval $[0,t]$ to deduce that  
\begin{equation}
\hat{\pi}_{t }^{\mu }-\hat{\pi}_{k\delta }^{\mu }
=\psi _{0 ,t }\left( \hat{\pi}_{0
}^{\mu }-\hat{x}_{0 }\right) +\psi _{0 ,t }\int_{0 }^{t}\psi _{0,s}^{-1}(P_{s}^{\mu }-P_{s})H^{\top }\mathrm{d}I_{s}  
+\psi _{0,t }\int_{0 }^{t}\psi _{0,s}^{-1}\mathrm{d}z_{s}^{\mu }.
\label{deltapi2'}
\end{equation}%
For a deterministic $\psi$, we can rewrite \eqref{deltapi2'} as 
\begin{equation}
\hat{\pi}_{t }^{\mu }-\hat{\pi}_{k\delta }^{\mu }
=\psi _{0 ,t }\left( \hat{\pi}_{0
}^{\mu }-\hat{x}_{0 }\right) +\int_{0 }^{t}\psi _{s,t}(P_{s}^{\mu }-P_{s})H^{\top }\mathrm{d}I_{s}  
+\int_{0 }^{t}\psi _{s,t}\mathrm{d}z_{s}^{\mu }.
\label{deltapi2''}
\end{equation}%
In this case, the bound of the stochastic terms in \eqref{deltapi2''} would follow directly from \eqref{exp} or, more precisely, from the deterministic exponential decay of $\psi_{s,t}$. However, in our case, $\psi$ is \emph{not} deterministic. Moreover the process $\psi_{s,t}$ is not adapted with respect to the filtration generated by the Brownian motion $I$.  To be more precise, it does not depend on $\{I_r, r\in [0,s]\}$ but on $\{I_r, r\in [0,t]\}$.      Hence, the stochastic integrals in \eqref{deltapi2''} do not make sense as standard It\^o integrals. One can interpret them using a more general definition of Skorohod integration, but the control of such resulting integrals is no longer immediately obvious. A bound using Malliavin calculus may be possible;  see the monograph of Nualart\cite{nualart} for details of the methodology.

In another approach, one could attempt to keep out of the stochastic integrals the part that is
not adapted, i.e. $\psi_{0,t}$, and only use \eqref{deltapi2'} but not \eqref{deltapi2''}. To do so, one would need to bound the exponential blow-up of the stochastic integrals in \eqref{deltapi2'} as $t\rightarrow \infty$ 
 over the entire positive half line $[0,\infty)$. Instead, we limit ourselves to apply Duhamel's priciple on intervals of the form $[k\delta, (k+1)\delta]$, where $\delta$ is small enough to be able to rely on the bound given by \eqref{psi-}.

Finally, note that, using again an argument based on the Gr\"onwall inequality, one obtains
that, for any $T>0$, there exists a constant $c_{T }$ such that the expected difference between the means of the FA process started from different initial distributions $\mu$ and $\pi_0$ satisties
\begin{equation}
\sup_{t\in [0,T]}\mathbb{E}\left[
\left\vert \hat{\pi}_{t}^{\mu }-\hat{\pi}_{t}\right\vert \right] \leq
c_{T }. \label{cd'}
\end{equation}%
Thus, we only need to check that there exists a time horizon $T$ for which
\begin{equation}
\sup_{t\in [T,\infty)}\mathbb{E}\left[
\left\vert \hat{\pi}_{t}^{\mu }-\hat{\pi}_{t}\right\vert \right] \leq
c_{T } \label{cd''}
\end{equation}%
The bound of the eighth moments of the covariance matrices $P_{t}^{\mu }$ and $P_{t}^{\pi _{0}}$,  respectively, in \eqref{controlpt}, as well as of the fourth moments of $\psi_{k\delta :s}^{-1}$, suffices for the result above, but it is perhaps not necessary. We leave it for future research to find the optimal bounds.       
\end{remark}

\begin{remark} \label{rem:linear} 
In the linear and Gaussian case, $\tilde f$ and $\tilde h$ are constant and $\lambda^\mu_s=0$. Therefore  
\[
z^{\mu}=\int_{0}^{t}(\tilde f-P_{s}^\mu H^{\top }\tilde h)\mathrm{d}s;
\] 
in other words, $z^\mu$ is a process with bounded variation. Moreover, $P_{s}^\mu$ is deterministic. In this case, all the technical difficulties described above vanish and Theorem \ref{nonlinear} holds true with a much simpler proof. In this case, though, a much stronger result holds, namely Theorem~\ref{linear}, which is proved in the next section of this appendix.          
\end{remark}

\begin{remark}\label{r3} 
Theorem \ref{nonlinear} gives us only an upper bound on the expected value of the distance between $\pi_{t}^\mu$ and  $\pi_t$. From it, we can derive the following alternative control: Let us denote by $T(t,R)$ the average time spent by the process $\pi^{\mu}_t $ outside the ball $B(\pi_t,R)$, that is \[
T(t,R)=\frac{1}{t}\int_0^tI_{\{W_{2}(\pi _{t}^{\mu },\pi_{t}^{0})>R\}}\mathrm{d}s.
\]  Then, for any $\varepsilon>0$, there exists a constant 
$R=R({\mu},\varepsilon)$ independent of $t$ such that   
\begin{equation*}
\sup_{t\ge 0}\mathbb E\left [T(t,R)\right]<\varepsilon.
\end{equation*}
We note that we cannot expect $T(t,R)$ to decrease to zero as $t$ tends to $\infty$, since $\hat\pi^{\mu}_t-\hat{x}_t $  can be viewed as a  $d$-dimensional Ornstein-Uhlenbeck process with random coefficients and perturbed  by the random residue process  $z$.      
\end{remark}

\begin{remark} The concept of stabilization described in Theorem \ref{nonlinear} is much weaker than the one described in many classical stabilization results, e.g. such as those covered by Budhiraja \cite{Budhiraja2011}, Picard\cite{picard1991} or Van Handel\cite{Handel.2009}. On the other hand, it has the advantage that the signal $X(t)$ is not required to be fully observable. In particular, the dimension $n$ of the observation $Y(t)$ and the dimension $d$ of the signal $X(t)$, in the notation of Sec.~\ref{sec:framework}, do not need to coincide. Heuristically we need to be able to observe all the "unstable" directions. See {\mg Section~D in this Appendix} 
for some simple illustrative examples in this direction. This feature of our results is very advantageous in practical applications as, in many situations, the dimension of the observation space is considerably smaller than that of the state space, $n \ll d$  \cite{Asch.Bocquet.2016,Bengtsson.ea.1981}.
\end{remark}
\begin{remark}One can deduce a criterion for the bound \eqref{controlpt} to hold. The derivation of such a criterion can be carried out in terms of the centralized third moments of $\pi_t^{\mu}$, along the lines of the arguments in Section 6.2 of Bain and Crisan \cite{bc}, and it is the subject of future research.
\end{remark}  

\subsection*{C. Assumptions and proof of Theorem~\ref{linear}} 

\label{sec:linear}  

The assumptions and the arguments in this section  are based on the framework of 
Ocone and Pardoux\cite{op} and, in part, on their results. In particular, we impose the 
following set of assumptions on the coefficients of the signal and observation equations 
\eqref{signal}+\eqref{observation} under which Theorem~\ref{linear} holds:

\begin{itemize}
\item We assume that the coefficients $f$ and $h$ are linear. 
In other words, we will assume that 
\begin{equation}
f=F{\mathcal{I}}+\tilde{f},\ \ h=H{\mathcal{I}}+\tilde{h},  \label{Fhsigmal}
\end{equation}%
where
\begin{enumerate}[(i)]
\item ${\mathcal{I}}:\mathbb{R}^{d}\rightarrow \mathbb{R}^{d}$ is the
identity function defined as ${\mathcal{I}}(x)=x$ for any $x\in \mathbb{R}%
^{d}$;
\item $F\in \mathbb{R}^{d\times d}$, $H\in \mathbb{R}^{d\times n}$ are
given matrices;
\item  $\tilde{f}$ and $\tilde{h}$ are $d$-dimensional and $n$%
-dimensional vectors, respectively; and
\item the function $\sigma $ is a constant $d\times d$-matrix.
\end{enumerate}

\item We assume that the measure $\mu $ has finite second
moments --- and thus it it belongs to the Wasserstein space, as discussed in 
{\mg Appendix I} --- and that it is absolutely continuous
with respect to $\pi _{0}$. We denote by $\theta _{\mu }$ the density of $%
\mu $ with respect to $\pi _{0}$ which is  integrable
with respect to $\pi _{0}.$ 
Using a standard probabilistic result (see, for
example, Problem 3.20 in Karatzas \& Shreve \cite{KaratzasShreve}), it follows that there
exists a random variable $\Upsilon _{\mu }\geq 0$ such that $%
\lim_{t\rightarrow \infty }\mathbb{E}\left[ \theta _{\mu }\left(
X_{0}\right) |\mathcal{Y}_{t}\right] =\Upsilon _{\mu }$. Moreover both $%
\mathbb{E}\left[ \theta _{\mu }\left( X_{0}\right) |\mathcal{Y}_{t}\right] $
and $\Upsilon _{\mu }$ \ have a pathwise representation in terms of the observation path $Y_{t}.$

We assume that, for the path for which we do the analysis, $\lim_{t\rightarrow \infty }%
\mathbb{E}\left[ \theta _{\mu }\left( X_{0}\right) |\mathcal{Y}_{t}\right]
=\Upsilon _{\mu }>0.$ Similarly, we will be assuming that both $\lim_{t\rightarrow
\infty }\mathbb{E}\left[ \theta _{\mu }\left( X_{0}\right) X_{0}|\mathcal{Y}%
_{t}\right] $ and $\lim_{t\rightarrow \infty }\mathbb{E}\left[ \theta _{\mu
}\left( X_{0}\right) X_{0}X_{0}^*|\mathcal{Y}_{t}\right] $ exist, 
based on the same standard probabilistic result.

\item We define $\pi _{t}^{0}:=N\left( \hat{x}_{t}^{0},P_{t}^{0}\right),$
where $\hat{x}$ satisfies the stochastic differential equation 
\begin{multline}
\mathrm{d}\hat{x}_{t}^{0} = \left( F\hat{x}_{t}^{0}+\tilde{f}\right) \,\mathrm{%
d}t+P_{t}H^{\top }(\mathrm{d}Y_{t}-(H\hat{x}_{t}^{0}+\tilde{h})\,\mathrm{d}%
t) \\
= (Q_{t}\hat{x}_{t}^{0}+\tilde{f}-P_{t}H^{\top }\tilde{h})\,\mathrm{d}%
t+P_{t}H_{t}^{\top }\mathrm{d}Y_{t},\quad \hat{x}_{0}=0
\label{eq:kalman:conditionalmean}
\end{multline}%
and $P^{0}$ satisfies the deterministic matrix Riccati equation 
\begin{multline}
\frac{\mathrm{d}P_{t}^{0}}{\mathrm{d}t} = \sigma _{t}\sigma _{t}^{\top
}+FP_{t}^{0}+P_{t}^{0}F^{\top }-P_{t}^{0}H^{\top}HP_{t}^{0}\\ 
= FQ_{t}+Q_{t}F^{\top }+\sigma _{t}\sigma _{t}^{\top
}+P_{t}H^{\top }HP_{t},  \label{eq:kalman:conditionalcovariance}
\end{multline}%
with $P_{0}=0\in \mathbb{R}^{d\times d}$ being the matrix with null entries.
We assume, furthermore, that there exists a unique solution $P_{\infty }\geq 0$ to
the algebraic Riccati equation 
\begin{equation}
\sigma \sigma ^{\top }+FP+PF^{\top }-PH^{\top }HP = 0,  \label{assym}
\end{equation}%
and that $Q_{t}:=F-P_{t}^{0}H^{\top }H$ is asymptotically stable.
Using this unique non-negative solution, we
define $Q_{\infty }:=F-P_{\infty }H^{\top }H$.
\cite{Handel.2009}\end{itemize}

Following from Remark 2.1 and Lemma 2.2 in Ocone and Pardoux\cite{op}, which
in turn uses the classical result of Theorem 3.7 in Kwakernaak and Sivan \cite%
{ks}, we make the following remark:
\begin{remark}\label{ptoinf} 
If $(F,H)$ is detectable and  $(F, \sigma_{0})$ is stabilizable, then \eqref{assym} has a unique non-negative definite solution $P_{\infty}$ and one has both $P_{t} \rightarrow P_{\infty}$ and $Q_{t} \rightarrow Q_{\infty}$  for any initial condition $P_{0}^{} \geq 0$, with $Q_t$ asymptotically stable. See Chapter 1 in Kwakernaak and Sivan \cite{ks} for the definition of detectability and stabilizability. In particular, the eigenvalues of $Q_\infty$ have negative real parts and  for any $0<a<\inf \left\{-\operatorname{Re} \lambda ; {\rm over} \,\, \lambda\right.$ {\rm being an eigenvalue of} $\left. Q_\infty \right\}$, there is a constant $K_{a}$ such that
$$
\left\|P_{t}^{}-P_{\infty}\right\| \leq K_{a} e^{-at}
$$
The last fact can be proved by observing that
\[
\frac{d}{d t}\left(P_{t}^{}-P_{\infty}\right)=\left[F-1 / 2\left(P_{t}^{}+P_{\infty}\right) H^{\top} H\right]\left(P_{t}^{}-P_{\infty}\right) 
+\left(P_{t}^{}-P_{\infty}\right)\left[F-1 / 2\left(P_{t}^{R}+P_{\infty}\right) H^{\top} H\right]^{\top}
\]
and carrying out an analysis similar to that in Theorem 2.3 of Ocone and Pardoux\cite{op}.
\end{remark}

We will show below that $\pi ^{\mu }$ gets asymptotically close to the
``reference process” $\pi ^{0}$, regardless of the initial condition, and
since this holds true for $\mu =\pi _{0}$ too, we immediately deduce that 
\begin{equation*}
\lim_{t\rightarrow \infty }W_{2}(\pi _{t}^{\mu },\pi _{t})=0.
\end{equation*}%
Note that $\pi $ coincides with $\pi ^{0}$ if the initial condition $\pi
_{0} $ is the Dirac delta distribution at $0$, that is, $\pi _{0}=\delta _{0}$. The
reference process $\pi ^{0}$ is convenient to work with: From Remark %
\ref{ptoinf} we deduce that the centered version of $\pi ^{0,-\hat{x}%
_{t}^{0}}:=N\left( 0,P_{t}^{0}\right) $ converges weakly, as well as in 
Wasserstein distance, to $\pi ^{\infty }=N\left( 0,P_{\infty }\right) $. Using
the equivalent definition of the Wasserstein distance 
\eqref{equiv}, it follows that 
\begin{equation*}
\lim_{t\rightarrow \infty }W_{2}(\pi _{t}^{\mu },\pi
_{t}^{0})=0\Longleftrightarrow \lim_{t\rightarrow \infty }W_{2}(\pi
_{t}^{\mu ,-\hat{x}_{t}^{0}},\pi _{t}^{0,-\hat{x}_{t}^{0}})=0%
\Longleftrightarrow \lim_{t\rightarrow \infty }W_{2}(\pi _{t}^{\mu ,-\hat{x}%
_{t}^{0}},\pi ^{\infty })=0.
\end{equation*}%

To prove the last limit, it suffices to show that $\pi _{t}^{\mu ,-\hat{x}%
_{t}^{0}}$ converges to $\pi ^{\infty }$ in the weak topology and the first
and second moments of $\pi _{t}^{\mu ,-\hat{x}_{t}^{0}}$ converge to the
first and second moments of $\pi ^{\infty }$, respectively. Equivalently, it
suffices to show that:
\begin{enumerate}[(a)]
    \item $\pi _{t}^{\mu ,-\hat{x}_{t}^{0}}$ converges to $\pi^{\infty }$ in the weak topology; 
    \item $\lim_{t\rightarrow \infty}| \hat{\pi}_{t}^{\mu }-\hat{x}_{t}^{0}| =0$ ,
   i.e., the distance betweeen the mean of $\pi _{t}^{\mu }$ and that of $\pi _{t}^{0}$
    tends to $0$; and 
    \item $\lim_{t\rightarrow \infty }| P_{\pi_{t}^{\mu }}-P_{t}^{0}| =0$, i.e., the distance 
    betweeen the covariance matrix of $\pi _{t}^{\mu }$ and that of $\pi _{t}^{0}$ tends to $0$.
\end{enumerate} 

Since the set of bounded uniformly continuous functions are convergence determining, to justify that 
$\pi _{t}^{\mu ,-\hat{x}_{t}^{0}}$ converges to $\pi ^{\infty }$ in the weak topology it suffices to
show that 
$\lim_{t\rightarrow \infty }|\pi _{t}^{\mu ,-%
\hat{x}_{t}^{0}}\left( \varphi \right) -\pi ^{\infty }\left( \varphi \right)
| =0$, for any $\varphi $ that is a bounded and uniformly continuous function. 
To recap, we have that $\lim_{t\rightarrow \infty}W_{2}(\pi _{t}^{\mu },\pi _{t}^{0})=0$, 
if and only if the following three
properties hold true:

\begin{itemize}
\item $\lim_{t\rightarrow \infty }| \hat{\pi}_{t}^{\mu }-\hat{x}%
_{t}^{0}| =0;$

\item $\lim_{t\rightarrow \infty }\left\vert P_{\pi _{t}^{\mu
}}-P_{t}^{0}\right\vert =0;$ and

\item $\lim_{t\rightarrow \infty }\left\vert \pi _{t}^{\mu }\left( \varphi
_{t}\right) -\pi _{t}^{0}\left( \varphi _{t}\right) \right\vert =0$ for any
bounded uniformly continuous function $\varphi $, where $\varphi
_{t}$ is the same function shifted by the mean $\hat{x}_{t}^{0},$ that is $%
\varphi _{t}\left( x\right) :=\varphi _{t}\left( x+\hat{x}_{t}^{0}\right) ,$ 
$x\in \mathbb{R}^{d}$. 
\end{itemize}
Note that the reference process $\pi ^{0}$ is not unique: we can replace
it by any other process $\tilde{\pi}^{0}$ with the property that $%
\lim_{t\rightarrow \infty }W_{2}(\tilde{\pi}_{t}^{0},\pi _{t}^{0})=0$.

The proof of Theorem \ref{linear} is largely based on the arguments in the
proof of Theorem 2.6 in Ocone and Pardoux\cite{op}.
Using the linearity of the equation~\eqref{signal}, we can 
decompose the signal as follows:
\begin{eqnarray*}
X_{t} &=&e^{Ft}X_{0}+\tilde{X}_{t}, \\
\tilde{X}_{t} &:&=\int_{0}^{t}(F\tilde{X}_{s}+f)\mathrm{d}%
s+\int_{0}^{t}\sigma \mathrm{d}V_{s}=\int_{0}^{t}e^{F(t-s)}(\tilde{f}%
\mathrm{d}s+\sigma \mathrm{d}V_{s}), \\
Y_{t} &=&\int_{0}^{t}(He^{Fs}X_{0}+H\tilde{X}_{s}+\tilde{h})\mathrm{d}%
s+W_{t},
\end{eqnarray*}%
and we introduce the new measure $\overline{\mathbb{P}}$ defined by 
\begin{equation*}
\left.\frac{\mathrm{d}\overline{\mathbb{P}}}{\mathrm{d}\mathbb{P}}\right |_{\mathcal{F}_t}=\exp \left[ \int_{0}^{t}-\left%
\langle He^{Fs}X_{0},\mathrm{d}W_{s}\right\rangle -\frac{1}{2}\int_{0}^{t}\left\vert
He^{Fs}X_{0}\right\vert ^{2}\mathrm{d}s\right],\ \ \ t\ge 0.
\end{equation*}%
By Girsanov's theorem, under this measure $\overline{\mathbb{P}},$ the process 
\begin{equation*}
\bar{W}_{t}:=\int_{0}^{t}He^{Bs}X_{0}\mathrm{d}s+W_{t},\quad t\leq T
\end{equation*}%
is a Brownian motion, and $X_{0}$ is independent of $\left( V_{t},\bar{W}%
_{t}\right) _{t\leq T}$.  Following the proof of Proposition 3.13  in Bain and Crisan\cite{bc}, the law of $X_{0}$ remains unchanged under $\bar {\mathbb P}$. Let 
\begin{equation*}
L_{t}:= \left.\frac{\mathrm{d}\mathbb{P}}{\mathrm{d}\overline{\mathbb{P}}}\right |_{\mathcal{F}_t}=\exp \left[ \int_{0}^{t}\left\langle He^{Bs}X_{0},\mathrm{d}\bar{W}%
_{s}\right\rangle -\frac{1}{2}\int_{0}^{t}\left\vert He^{Bs}X_{0}\right\vert
^{2}\mathrm{d}s\right] =\exp \left( \left\langle X_{0},W^{\circ }\right\rangle -\frac{%
1}{2}\left\langle M_{t}X_{0},X_{0}\right\rangle \right) ,
\end{equation*}%
where 
\begin{equation*}
\left( W_{t}^{o},M_{t}\right) =\left( e^{F^*t}\int_{0}^{t}H^*\mathrm{d}\bar{W}%
_{s},\int_{0}^{t}e^{F^*s}H^*He^{Fs}\mathrm{d}s\right).
\end{equation*}%
Let $%
\psi :$ $\mathbb{R}^{d}\times \mathbb{R}^{d}\rightarrow \mathbb{R}$ be a function that is
integrable with respect to the joint law of $\left( X_{0},\tilde{X}%
_{t}\right)$. Then 
\begin{equation*}
\mathbb{E}\left[ \psi \left( X_{0},\tilde{X}_{t}\right) |\mathcal{Y}_{t}%
\right] =\frac{\bar{\mathbb{E}}\left[ \psi \left( X_{0},\tilde{X}_{t}\right)
L_{t}|\mathcal{Y}_{t}\right] }{\bar{\mathbb{E}}\left[ L_{t}|\mathcal{Y}_{t}%
\right] }.
\end{equation*}%
From the above formula for $\tilde{X}_{t}$ and the fact that $X_{0}$ and $(%
\bar{W},Y)$ are $\overline{\mathbb{P}}$-independent, we deduce that

\begin{equation}
\mathbb{E}\left[ \psi \left( X_{0},\tilde{X}_{t}\right) |\mathcal{Y}_{t}%
\right] =\frac{\int_{\mathbb{R}^{d}}e^{-1/2\left( M_{t}x,x\right) }\bar{%
\mathbb{E}}\left[ \psi \left( x,\tilde{X}_{t}\right) \exp \left(
\left\langle x,W^{\circ }\right\rangle \right) |\mathcal{Y}_{t}\right] \pi
_{0}(\mathrm{d}x)}{\int_{\mathbb{R}^{d}}e^{-1/2\left( M_{t}x,x\right) }\bar{\mathbb{E}%
}\left[ \exp \left( \left\langle x,W^{\circ }\right\rangle \right) |\mathcal{%
Y}_{t}\right] \pi _{0}(\mathrm{d}x)}.  \label{pair}
\end{equation}%
In particular, for any function $\varphi :\mathbb{R}^{d}\rightarrow \mathbb{R%
}$ such that the random variable $\varphi \left( X_{t}\right) =\varphi
\left( e^{Ft}X_{0}+\tilde{X}_{t}\right) $ is integrable, we deduce from (\ref%
{pair}) that%
\begin{eqnarray*}
\mathbb{E}\left[ \varphi \left( X_{t}\right) |\mathcal{Y}_{t}\right] &=&%
\frac{\int_{\mathbb{R}^{d}}e^{-1/2\left( M_{t}x,x\right) }\bar{\mathbb{E}}%
\left[ \varphi _{x}^{1}\left( \tilde{X}_{t},W^{\circ }\right) |\mathcal{Y}%
_{t}\right] \pi _{0}(\mathrm{d}x)}{\int_{\mathbb{R}^{d}}e^{-1/2\left( M_{t}x,x\right)
}\bar{\mathbb{E}}\left[ \varphi _{x}^{2}\left( \tilde{X}_{t},W^{\circ
}\right) |\mathcal{Y}_{t}\right] \pi _{0}(\mathrm{d}x)}; \\
\varphi _{x}^{1}\left( a,b\right) &=&\varphi \left( e^{Ft}x+a\right) \exp
\left\langle x,b\right\rangle ,~~a\in \mathbb{R}^{d}, \, b\in \mathbb{R}^{n}; \,\, {\rm and} \\
\varphi _{x}^{2}\left( a,b\right) &=&\exp \left\langle x,b\right\rangle,
~~~~a\in \mathbb{R}^{d}, \, b\in \mathbb{R}^{n}.
\end{eqnarray*}%

The pair $\left( \tilde{X},W^{\circ }\right) $ satisfies a linear system of
stochastic differential equations driven by the $\left( d+m\right) $-Brownian motion $\left( V,W^{\circ
}\right) $ with initial condition $\left( 0,0\right) $. Moreover, under $%
\overline{\mathbb{P}},$ the process $Y$ satisfies 
\begin{equation}\label{newobs}
Y_{t}=\int_{0}^{t}(H\tilde{X}_{s}+\tilde{h})\mathrm{d}s+\bar{W}_{t}.
\end{equation}%
It follows that we can express $\bar{\mathbb{E}}\left[ \varphi
_{x}^{i}\left( \tilde{X}_{t},W_{t}^{\circ }\right) |\mathcal{Y}_{t}\right]
,i=1,2,$ as integrals with respect to a Gaussian distribution $\eta $  
 with mean vector  $\left( \hat{x}_{t}^{0},\hat{x}%
_{t}^{W}\right)^*$ and covariance matrix%
\begin{equation*}
C_{t}=\left( 
\begin{array}{cc}
P_{t}^{0} & S_{t} \\ 
S_{t} & R_{t}%
\end{array}%
\right) 
\end{equation*}%
that satisfy the equations satisfied by a Kalman-Bucy filter with the signal equations identical to those for the pair 
$\left( \tilde{X}_{t},W_{t}^{\circ }\right)$ and the observation equation identical to (\ref{newobs}). 

More precisely,   $\hat{x}_{t}^{0}$ is the solution of equation (\ref%
{eq:kalman:conditionalmean}), $P_{t}^{0}$ is the solution of the equation (\ref%
{eq:kalman:conditionalcovariance}), $\hat{x}_{t}^{W}$ is the solution of the
equation 
\begin{equation*}
\mathrm{d}\hat{x}_{t}^{W}=\left( e^{Ft}+S_{t}\right) ^{\ast }H^{\ast }\left(
\mathrm{d}Y_{t}-\left( H\hat{x}_{t}^{0}+\tilde{h}\right) \mathrm{d}t\right) ,\quad \hat{x}%
_{0}^{W}=0
\end{equation*}%
and the pair $\left( S_{t},Q_{t}\right) $
solve the equation s
\begin{eqnarray*}
\frac{\mathrm{d}S_{t}}{\mathrm{d}t} &=&FS_{t}-P_{t}^{0}H^{\ast }H\left( e^{Ft}+S_{t}\right)
,\quad S_{0}=0, \\
\frac{\mathrm{d}R_{t}}{\mathrm{d}t} &=&-e^{F^{\ast }t}H^{\ast }HS_{t}-S_{t}^{\ast }H^{\ast
}He^{Ft}-S_{t}^{\ast }H^{\ast }HS_{t},\quad R_{0}=0;
\end{eqnarray*}%
see formulae (27--33) in Ocone and Pardoux\cite{op} for details. 
One can check that 
\begin{equation*}
\frac{\mathrm{d}S^\circ_t}{\mathrm{d}t} =\left( F-P_{t}^{0}H^{\ast }H\right)
S^\circ_t=Q_{t}S^\circ_t,
\end{equation*}%
where $S^\circ_t:=e^{Ft}+S_{t}$ and, since $Q_{t}$ is asymptotically stable, we conclude that 
\begin{equation}
\lim_{t\rightarrow \infty }\left\vert \left\vert S^\circ_t\right\vert
\right\vert =0.  \label{eq:sumas}
\end{equation}%
In fact, similar to Remark \ref{ptoinf}, the convergence in (\ref{eq:sumas}) has an exponential decay rate to 0. 

We  deduce that 
\begin{eqnarray}
\bar{\mathbb{E}}\left[ \varphi _{x}^{2}\left( \tilde{X}_{t},W^{\circ
}\right) |\mathcal{Y}_{t}\right] &=&\bar{\mathbb{E}}\left[ \exp \left\langle
x,W^{\circ }\right\rangle |\mathcal{Y}_{t}\right] = e^{\left\langle x,\hat{x}%
_{t}^{W}\right\rangle +1/2\left( Q_{t}x,x\right) } ; \label{help} \\
\bar{\mathbb{E}}\left[ \varphi _{x}^{1}\left( \tilde{X}_{t},W^{\circ
}\right) |\mathcal{Y}_{t}\right] &=&\bar{\mathbb{E}}\left[ \varphi \left(
e^{Ft}x+\tilde{X}_{t}\right) \exp \left\langle x,W^{\circ }\right\rangle |%
\mathcal{Y}_{t}\right]  \notag \\
&=&e^{\left\langle x,\hat{x}_{t}^{W}\right\rangle +1/2\left( Q_{t}x,x\right)
}\int \varphi \left( S^\circ_tx+a\right) \pi
_{t}^{0}\left( \mathrm{d}a\right) .  \notag
\end{eqnarray}%
It follows that 
\begin{equation*}
\pi _{t}\left( \varphi \right) =\mathbb{E}\left[ \varphi \left( X_{t}\right)
/\mathcal{Y}_{t}\right] =\frac{\int_{\mathbb{R}^{d}}\Xi _{t}\left( x\right)
\int \varphi \left( S^\circ_tx+a\right) \pi _{t}^{0}\left(
\mathrm{d}a\right) \pi _{0}(\mathrm{d}x)}{\int_{\mathbb{R}^{d}}\Xi _{t}\left( x\right) \pi
_{0}(\mathrm{d}x)},
\end{equation*}%
where $\Xi _{t}\left( x\right) :=e^{-1/2\left( M_{t}x,x\right) +\left\langle
x,\hat{x}_{t}^{W}\right\rangle +1/2\left( Q_{t}x,x\right) }.$ Moreover, by
choosing a function $\psi :$ $\mathbb{R}^{d}\times \mathbb{R}^{d}\rightarrow 
\mathbb{R}$ in (\ref{pair}) that is independent in the second component, $%
\psi \left( x,y\right) =\varphi \left( x\right) $, $\left( x,y\right) \in 
\mathbb{R}^{d}\times \mathbb{R}^{d}$, we obtain from (\ref{help}) that 
\begin{equation*}
\mathbb{E}\left[ \varphi \left( X_{0}\right) |\mathcal{Y}_{t}\right] =\frac{%
\int_{\mathbb{R}^{d}}\Xi _{t}\left( x\right) \varphi \left( x\right) \pi
_{0}(\mathrm{d}x)}{\int_{\mathbb{R}^{d}}\Xi _{t}\left( x\right) \pi _{0}(\mathrm{d}x)}=\frac{%
\pi _{0}\left( \Xi _{t}\varphi \right) }{\pi _{0}\left( \Xi _{t}\right) }.
\end{equation*}

We identify the FA process $\lambda :[0,\infty )\ \times 
\mathcal{P}(\mathbb{R}^{d})\times \Omega \rightarrow \mathcal{P}(\mathbb{R}%
^{d})$ by the formula%
\begin{eqnarray}
\pi _{t}^{\mu }\left( \varphi \right) &=&\lambda (t,\omega )\mu \left(
\varphi \right)  \notag \\
&:&=\frac{\int_{\mathbb{R}^{d}}\Xi _{t}\left( x\right) \int \varphi \left(
S^\circ_t x+a\right) \pi _{t}^{0}\left( \mathrm{d}a\right) \mu (\mathrm{d}x)}{%
\int_{\mathbb{R}^{d}}\Xi _{t}\left( x\right) \mu (\mathrm{d}x)}  \notag \\
&=&\frac{1}{C^{\pi _{t}^{\mu }}}\int_{\mathbb{R}^{d}}\Xi _{t}\left( x\right)
\int \varphi \left( S^\circ_tx+a\right) \pi _{t}^{0}\left(
\mathrm{d}a\right) \theta _{\mu }\left( x\right) \pi _{0}(\mathrm{d}x),  \label{FAlin}
\end{eqnarray}%
where $C^{\pi _{t}^{\mu }}$ is the normalization constant 
\begin{equation*}
C^{\pi _{t}^{\mu }}:=\int_{\mathbb{R}^{d}}\Xi _{t}\left( x\right) \theta
_{\mu }\left( x\right) \pi _{0}(\mathrm{d}x)=\pi _{0}\left( \Xi _{t}\theta _{\mu
}\right) =\mathbb{E}\left[ \theta
_{\mu } \left( X_{0}\right) |\mathcal{Y}_{t}%
\right] \pi _{0}\left( \Xi _{t}\right) .
\end{equation*}%
To obtain the representation (\ref{FAlin}), we used the fact that $\mu $ is
absolutely continuous with respect to $\pi _{0}$ and that $\theta _{\mu }$ is
the density of $\mu $ with respect to $\pi _{0}$.

\paragraph*{Convergence of the first moments.} Observe that 
\begin{eqnarray}
\hat{\pi}_{t}^{\mu } &=&\frac{1}{C^{\pi _{t}^{\mu }}}\int_{\mathbb{R}%
^{d}}\Xi \left( x\right) \int (S^\circ_tx+a)\pi
_{t}^{0}\left( \mathrm{d}a\right) \theta _{\mu }\left( x\right) \pi _{0}(\mathrm{d}x)  \notag
\\
&=&\hat{x}_{t}^{0}+\frac{S^\circ_t \int_{\mathbb{R}%
^{d}}x\Xi _{t}\left( x\right) \theta _{\mu }\left( x\right) \pi _{0}(\mathrm{d}x)}{%
\int_{\mathbb{R}^{d}}\Xi _{t}\left( x\right) \theta _{\mu }\left( x\right)
\pi _{0}(\mathrm{d}x)}=\hat{x}_{t}^{0}+\frac{S^\circ_t \left( 
\mathbb{E}\left[ \theta _{\mu }\left( X_{0}\right) X_{0}|\mathcal{Y}_{t}%
\right] \pi _{0}\left( \Xi _{t}\right) \right) }{\left( \mathbb{E}\left[
\theta _{\mu }\left( X_{0}\right) |\mathcal{Y}_{t}\right] \pi _{0}\left( \Xi
_{t}\right) \right) }.  \label{meanhat}
\end{eqnarray}%
It follows that 
\begin{equation*}
\left\vert \hat{\pi}_{t}^{\mu }-\hat{x}_{t}^{0}\right\vert \leq \frac{\left\vert
\left\vert S^\circ_t\right\vert \right\vert \left\vert \left\vert \mathbb{%
E}\left[ \theta _{\mu }\left( X_{0}\right) X_{0}|\mathcal{Y}_{t}\right]
\right\vert \right\vert }{\mathbb{E}\left[ \theta _{\mu }\left( X_{0}\right)
|\mathcal{Y}_{t}\right] }.
\end{equation*}%
Since the processes $t\rightarrow $ $\mathbb{E}\left[ \theta _{\mu }\left(
X_{0}\right) X_{0}|\mathcal{Y}_{t}\right] $ and $t\rightarrow $ $\mathbb{E}%
\left[ \theta _{\mu }\left( X_{0}\right) |\mathcal{Y}_{t}\right] $ converge
and the second limit is positive, we have, by using (\ref{eq:sumas}), 
\begin{equation}\label{com} 
\lim_{t\rightarrow \infty }\left\vert \hat{\pi}_{t}^{\mu }-\hat{x}%
_{t}^{0}\right\vert =\lim_{t\rightarrow \infty }\left\vert \left\vert
S^\circ_t\right\vert \right\vert \frac{\lim_{t\rightarrow \infty
}\left\vert \left\vert \mathbb{E}\left[ \theta _{\mu }\left( X_{0}\right)
X_{0}|\mathcal{Y}_{t}\right] \right\vert \right\vert }{\lim_{t\rightarrow
\infty }\left\vert \left\vert \mathbb{E}\left[ \theta _{\mu }\left(
X_{0}\right)|\mathcal{Y}_{t}\right] \right\vert \right\vert }=0.
\end{equation}%
Again, similar to Remark \ref{ptoinf}, the convergence in (\ref{com}) is exponentially fast. Since $\lim_{t\rightarrow \infty }e^{-\bar{\epsilon}t}\left\vert \pi
_{t}^{\mu }\right\vert =\lim_{t\rightarrow \infty }e^{-\bar{\epsilon}%
t}\left\vert \hat{x}_{t}^{0}\right\vert =0$, for any $\bar{\epsilon}>0$, it
follows that, for any $i,j=1,...,d$, 
\begin{equation*}
\lim_{t\rightarrow \infty }\left\vert \left( \pi _{t}^{\mu }\right)
^{i}\left( \pi _{t}^{\mu }\right) ^{j}-\left( \hat{x}_{t}^{0}\right)
^{i}\left( \hat{x}_{t}^{0}\right) ^{j}\right\vert =0.
\end{equation*}%
The last limit is used in the proof of the convergence of the covariance matrix below.

\paragraph*{Convergence of the covariance matrix.} Choose%
\begin{equation*}
\varphi ^{ij}:\mathbb{R}^{d}\rightarrow \mathbb{R},~~~\varphi ^{ij}\left(
x^{1},..,x^{d}\right) :=x^{i}x^{j}-\left( \pi _{t}^{\mu }\right) ^{i}\left(
\pi _{t}^{\mu }\right) ^{j}
\end{equation*}%
We have that 
\begin{eqnarray*}
( P_{\pi _{t}^{\mu }}) ^{ij}-\left( P_{t}^{0}\right) ^{ij} &=&\pi
_{t}^{\mu }\left( \varphi ^{ij}\right) -\left( P_{t}^{0}\right) ^{ij} \\
&=&\frac{1}{C^{\pi _{t}^{\mu }}}\int_{\mathbb{R}^{d}}\Xi \left( x\right)
\int \varphi ^{ij}\left( \left( S^\circ_t\right) x+a\right) \pi
_{t}^{0}\left( da\right) \theta _{\mu }\left( x\right) \pi _{0}(dx)-\left(
P_{t}^{0}\right) ^{ij} \\
&=&\sum_{k,l}\frac{\left( S^\circ_t\right) ^{il}\left(
S^\circ_t\right) ^{jk}\left( \mathbb{E}_{\mu }\left[ \theta _{\mu }\left(
X_{0}\right) X_{0}^{l}X_{0}^{k}|\mathcal{Y}_{t}\right] \right) }{\mathbb{E}%
\left[ \theta _{\mu }\left( X_{0}\right) |\mathcal{Y}_{t}\right] } \\
&&+\sum_{l}\frac{\left( S^\circ_t\right) ^{il}\left( \mathbb{E}_{\mu }%
\left[ \theta _{\mu }\left( X_{0}\right) X_{0}^{l}|\mathcal{Y}_{t}\right]
\right) }{\mathbb{E}\left[ \theta _{\mu }\left( X_{0}\right) |\mathcal{Y}_{t}%
\right] }\left( \hat{x}_{t}^{0}\right) ^{j} \\
&&+\sum_{l}\frac{\left( S^\circ_t\right) ^{jl}\left( \mathbb{E}_{\mu }%
\left[ \theta _{\mu }\left( X_{0}\right) X_{0}^{l}|\mathcal{Y}_{t}\right]
\right) }{\mathbb{E}\left[ \theta _{\mu }\left( X_{0}\right) |\mathcal{Y}_{t}%
\right] }\left( \hat{x}_{t}^{0}\right) ^{i} \\
&&+\left( \hat{x}_{t}^{0}\right) ^{i}\left( \hat{x}_{t}^{0}\right)
^{j}-\left( \pi _{t}^{\mu }\right) ^{i}\left( \pi _{t}^{\mu }\right) ^{j},
\end{eqnarray*}%
which gives the required convergence to 0.

\paragraph*{Convergence for bounded uniformly continuous test functions.} We
have that
\begin{equation*}
\left\vert \pi _{t}^{\mu }\left( \varphi _{t}\right) -\pi _{t}^{0}\left(
\varphi _{t}\right) \right\vert \leq \frac{1}{C^{\pi _{t}^{\mu }}}\int_{%
\mathbb{R}^{d}}\Xi \left( x\right) \int \left\vert \varphi _{t}\left( \left(
S^\circ_t\right) x+a\right) -\varphi _{t}\left( a\right) \right\vert \pi
_{t}^{0}\left( da\right) \theta _{\mu }\left( x\right) \pi _{0}(dx).
\end{equation*}%
Decomposing next the integral in the numerator into the sum of the
integral over the region $\left\vert ( S^\circ_t) x\right\vert
<\delta $ and the integral over the region $\left\vert \left(
S^\circ_t\right) x\right\vert \geq \delta $ yields

\begin{eqnarray*}
\left\vert \pi _{t}^{\mu }\left( \varphi _{t}\right) -\pi _{t}^{0}\left(
\varphi _{t}\right) \right\vert  &\leq &\sup_{\left\vert y-y^{\prime
}\right\vert <\delta }\left\vert \varphi (y)-\varphi \left( y^{\prime
}\right) \right\vert +2\Vert \varphi \Vert _{\infty }\frac{\left( \mathbb{E}%
_{\mu }\left[ \theta _{\mu }\left( X_{0}\right) \mathbf{1}_{\left\{
\left\vert \left( S^\circ_t\right) X_{0}\right\vert \geq \delta \right\}
}|\mathcal{Y}_{t}\right] \right) }{\mathbb{E}\left[ \theta _{\mu }\left(
X_{0}\right) |\mathcal{Y}_{t}\right] } \\
&\leq &\sup_{\left\vert y-y^{\prime }\right\vert <\delta }\left\vert \varphi
(y)-\varphi \left( y^{\prime }\right) \right\vert +2\frac{\Vert \varphi
\Vert _{\infty }}{\delta ^{2}}\frac{\left( S^\circ_t\right) \mathbb{E}%
_{\mu }\left[ \theta _{\mu }\left( X_{0}\right) \left\vert X_{0}\right\vert
^{2}|\mathcal{Y}_{t}\right] }{\mathbb{E}\left[ \theta _{\mu }\left(
X_{0}\right) |\mathcal{Y}_{t}\right] }.
\end{eqnarray*}

As above, it follows that 
\begin{equation*}
\limsup_{t\rightarrow \infty }\left\vert \pi _{t}^{\mu }\left( \varphi
_{t}\right) -\pi _{t}^{0}\left( \varphi _{t}\right) \right\vert \leq
\sup_{\left\vert y-y^{\prime }\right\vert <\delta }\left\vert \varphi
(y)-\varphi \left( y^{\prime }\right) \right\vert 
\end{equation*}
and then using the uniform continuity of $\varphi $, we deduce the result and
thus complete the proof. \qed

\subsection*{D. A simple example}

\label{simpleexamples}

In this subsection, we provide an illustrative example of unstable dynamics and partial observations that still yield convergence of the FA process; see also refs.~\cite{Carrassi.ea.2008, Carrassi.ea.2007}. We choose a two-dimensional linear signal $\left( X^{1},X^{2}\right) $ with $%
\lambda _{1}<0$ and $\lambda _{2}>0$: 
\begin{equation*}
X_{t}^{i} =x_{0}^{i}+\int_{0}^{t}\lambda _{i}X_{s}^{i}\mathrm{d}%
s+\int_{0}^{t}\sigma ^{i}\mathrm{d}V_{s}^{i} ,\ \ i=1,2, 
\end{equation*}
so that the law of $\left( X^{1},X^{2}\right) $ is not stable. More precisely,
if we choose two systems $\left( X^{1},X^{2}\right) $, $\left( \tilde{X}^{1},%
\tilde{X}^{2}\right) $ starting from $\left( x_{0}^{1},x_{0}^{2}\right) $ and $%
\left( \tilde{x}_{0}^{1},\tilde{x}_{0}^{2}\right) $, respectively, then
their corresponding expected values drift away from each other. In
particular, 
\begin{equation*}
\lim_{t\rightarrow \infty }|E\left[ X_{t}^{2}\right] -E\left[ X_{t}^{2}%
\right] |=\lim_{t\rightarrow \infty }e^{\lambda _{2}t}\left\vert x_{0}^{2}-%
\tilde{x}_{0}^{2}\right\vert =\infty .
\end{equation*}
As a result, the Wasserstein distance $W_2$ between $p_t$ and $\tilde p_t$ tends
to $\infty$. 

In fact, one can easily show that the first component is stable, whilst
the second one is unstable. However, we can stabilize the system by observing 
\emph{only} the second component, namely the unstable one. For example, we can
choose a one-dimensional observation process of the form%
\begin{equation*}
dY_{t}=hX_{s}^{2}\,\mathrm{d}t+dW_{t},
\end{equation*}%
which will guarantee that
\begin{equation*}
\lim_{t\rightarrow \infty }d_{W}\left( \pi _{t},\tilde{\pi}_{t}\right) = 0.
\end{equation*}%
One can justify this either by checking that the assumptions of Theorem \ref%
{linear} are satisfied or by explicit calculation. Using the second approach, one notices that the probability
measures $\pi _{t}$ and $\tilde{\pi}_{t}$ are both Gaussian and 
shows that (i) the distance between the coresponding means converges to $0$ and
(ii) that the covariance matrices of $\pi _{t}$ and $\tilde{\pi}_{t}$
coincide. The latter two matrices are given by diag$\left( q_{t}^{11},q_{t}^{22}\right) $. They
are diagonal and 
\begin{eqnarray*}
\frac{\mathrm{d}}{\mathrm{d}t}\left( 
\begin{array}{cc}
q_{t}^{11} & 0 \\ 
0 & q_{t}^{22}%
\end{array}%
\right) &=&\left( 
\begin{array}{cc}
1 & 0 \\ 
0 & 1%
\end{array}%
\right) +\left( 
\begin{array}{cc}
\lambda _{1} & 0 \\ 
0 & \lambda _{2}%
\end{array}%
\right) \left( 
\begin{array}{cc}
q_{t}^{11} & 0 \\ 
0 & q_{t}^{22}%
\end{array}%
\right) +\left( 
\begin{array}{cc}
q_{t}^{11} & 0 \\ 
0 & q_{t}^{22}%
\end{array}%
\right) \left( 
\begin{array}{cc}
\lambda _{1} & 0 \\ 
0 & \lambda _{2}%
\end{array}%
\right) \\
&&-\left( 
\begin{array}{cc}
q_{t}^{11} & 0 \\ 
0 & q_{t}^{22}%
\end{array}%
\right) \left( 
\begin{array}{cc}
0 & 0 \\ 
0 & h^{2}%
\end{array}%
\right) \left( 
\begin{array}{cc}
q_{t}^{11} & 0 \\ 
0 & q_{t}^{22}%
\end{array}%
\right), \\
\frac{\mathrm{d}}{\mathrm{d}t}\left( 
\begin{array}{cc}
q_{t}^{11} & 0 \\ 
0 & q_{t}^{22}%
\end{array}%
\right) &=&\left( 
\begin{array}{cc}
1+\lambda _{1}p_{t}^{22} & 0 \\ 
0 & 1+2\lambda _{2}p_{t}^{11}-h^{2}\left( p_{t}^{11}\right) ^{2}%
\end{array}%
\right),
\end{eqnarray*}%
which implies that 
\begin{equation*}
\lim_{t\rightarrow \infty }\left( 
\begin{array}{c}
p_{t}^{11} \\ 
p_{t}^{22}%
\end{array}%
\right) = \left( 
\begin{array}{c}
- {\displaystyle \frac{1}{2\lambda _{1}}} \\ 
{\displaystyle \frac{ \sqrt{\left( \lambda _{2}\right) ^{2}+h^{2}}-\lambda _{2}}{h^{2}}}%
\end{array}%
\right)
\end{equation*}%
and the limit is valid independent of the initial condition. This result
together with the Gaussianity property of the processes $\pi _{t}$ and $%
\tilde{\pi}_{t}$ implies the convergence in Wasserstein distance.

In the reverse situation, in which we would observe the stable component and
not the unstable one, the posterior distribution would not be stable.

In the above example, one can still apply existing results to the
one-dimensional unstable system and, after coupling it to the unobserved
component, obtain the stability of the pair. This works because the
observation process depends only on the observed component. In multidimensional,
nonlinear applications with very large dimensions $d$ and $n$, though,
it would be quite difficult to do this as it may not be possible to identify 
the unstable components a priori. 

The result presented in this paper does not assume that one needs to identify 
a system's unstable components and observe those. For example, if 
the observation process has the form%
\begin{equation*}
dY_{t}=h\left( X_{s}^{1}+X_{s}^{2}\right) \mathrm{d}t+\mathrm{d}W_{t},
\end{equation*}%
the FA process will still stabilize the system. More precisely, it is still the case 
that both $\pi _{t}$ and $\tilde{\pi}%
_{t}$ are Gaussian and that
\begin{equation*}
\lim_{t\rightarrow \infty }d_{W}\left( \pi _{t},\tilde{\pi}_{t}\right) = 0.
\end{equation*}

\subsection*{E. The forecast--assimilation (FA) process as a random dynamical
system (RDS)}

\label{sec:rds}

In this subsection, we give a brief justification of the fact that the solution
of \eqref{ks0} can be recast as an RDS evolving in the space of
probability measures $\mathcal{P}(\mathbb{R}^{d})$. This result is
related to similar  approaches in the nonlinear filtering literature; see,
for instance, Budhiraja\cite{Budhiraja2011}, Kunita\cite{Kunita.1986}, Picard\cite{picard1991}, 
Ocone and Pardoux\cite{op}, etc. The cornerstone of the argument is that the solution of 
equation (\ref{ks0}) can be expressed as 
\begin{equation}
\pi _{t}^{\mu }=\lambda (t,\omega )\mu =\frac{1}{C^{\pi _{0}}}\rho
_{0:t}(\mu ),  \label{kallstrie}
\end{equation}%
where $\rho _{t_{0},t_{1}}(\mu )$ is a two-parameter measure-valued process
defined as 
\begin{equation}
\rho _{t_{0},t_{1}}(\mu )\left( \varphi \right) =\int_{\mathbb{R}^{d}}%
\mathbb{\tilde{E}}_{t_{0},x}[\tilde{Z}_{t_{0},t_{1}}\varphi (\tilde{X}%
_{t_{0},t_{1}})\mid \mathcal{Y}]\mu \left( dx\right) ,  \label{ztilde}
\end{equation}%
where $\varphi $ is an arbitrary bounded Borel measurable map and $C^{\mu }$ is
the normalization constant $C^{\mu }=\rho _{t_{0},t_{1}}(\mu )(1)$. 
In (\ref{ztilde}), we have the following definitions:
\begin{enumerate}[(i)]
\item the process $\tilde{Z}_{t_{0},t_{1}}=\{\tilde{Z}_{t_{0},t_{1}},\ t\geq
0\}$ is defined by
\begin{equation}
\tilde{Z}_{t_{0},t_{1}}=\exp \left( \sum_{i=1}^{m}\int_{t_{0}}^{t_{1}}h^{i}(%
\tilde{X}_{t_{0},s})\,\mathrm{d}Y_{s}^{i}-\frac{1}{2}\sum_{i=1}^{m}%
\int_{t_{0}}^{t_{1}}h^{i}(\tilde{X}_{t_{0},s})^{2}\,\mathrm{d}s\right) ;
\label{ztilde2}
\end{equation}

\item the process $\tilde{X}_{t_{0},(\cdot) }=\{\tilde{X}_{t_{0},t_{1}},\
t_{0},t_{1}\geq 0\}$ is a stochastic process independent of $Y$ satisfying
the signal equation (\ref{signal}) on $[t_{0},\infty )$; and

\item $\mathbb{\tilde{E}}_{t_{0},x}$ is the expectation with respect to a
probability measure $\mathbb{P}_{t_{0},x}$ under which $Y$ is a Brownian
motion independent of $\tilde{X}$ and $\tilde{X}%
_{t_{0},t_{0}}\equiv x$.
\end{enumerate}

The independence of $\tilde{X}$ and $Y$ under $\mathbb{P}_{t_{0},x}$ in
formula (\ref{kallstrie}) enables us to show that $\rho _{t_{0},t_{1}}$ is an
RDS, which will immediately imply that $\pi _{t_{0},t_{1}}
$ is one, too. To justify this we introduce $\Theta ^{y_{(\cdot) }}(t_{0},t_{1})
$ to be the following two-parameter family of random variables 
\begin{equation}
\Theta ^{y_{(\cdot) }}(t_{0},t_{1})\triangleq \exp \!\left( h(\tilde{X}%
_{t_{0},t_{1}})^{\top }y_{t_{1}}-h(\tilde{X}_{t_{0},t_{0}})^{\top
}y_{t_{0}}+I_{t_{0},t_{1}}^{y_{(\cdot) }}-\frac{1}{2}\sum_{i=1}^{m}%
\int_{t_{0}}^{t_{1}}h^{i}(\tilde{X}_{t_{0},s})^{2}\,\mathrm{d}s\right) ,
\label{eq:robust:babar}
\end{equation}%
where $I_{t_{0},t_{1}}^{y_{(\cdot) }}$, is a version of the stochastic
integral $\int_{t_{0}}^{t_{1}}y_{s}^{\top }\,\mathrm{d}h(\tilde{X}_{t_{0},s})
$ and $y_{(\cdot) }$ is a continuous path, $y_{(\cdot) }\in C_{\mathbb{R}%
^{m}}[0,\infty ).$ The argument of the exponent in the definition of $\Theta
^{y_{(\cdot) }}(t_{0},t_{1})$ is recognizable as a formal integration by
parts of the argument of the exponential in (\ref{ztilde}). 


Let $\rho _{t_{0},t_{1}}^{y_{(\cdot) }}\left( \mu \right) $ and $\pi
_{t_{0},t_{1}}^{y_{(\cdot) }}(\mu )$ be the following two-parameter measure
valued processes, 
\begin{subequations}  \label{eq:robust}
        \begin{align}
\rho _{t_{0},t_{1}}^{y_{(\cdot) }}\left( \mu \right) (\varphi ) &= \int_{%
\mathbb{R}^{d}}\mathbb{\tilde{E}}_{t_{0},x}[\varphi (\tilde{X}%
_{t_{0},t_{1}})\Theta ^{y_{(\cdot) }}(t_{0},t_{1})\mid \mathcal{Y}]\mu \left(
dx\right),   \label{eq:robust:fdefn} \\
\pi _{t_{0},t_{1}}^{y_{(\cdot) }}(\mu )(\varphi ) &= \frac{\rho
_{t_{0},t_{1}}^{y_{(\cdot) }}\left( \mu \right) (\varphi )}{\rho
_{t_{0},t_{1}}^{y_{(\cdot) }}\left( \mu \right) (1)}.  \label{eq:robust:pidefn}
        \end{align}
\end{subequations}%
Then $\rho _{t_{0},t_{1}}^{Y_{(\cdot) }}(\mu )$ and $\pi
_{t_{0},t_{1}}^{Y_{(\cdot) }}(\mu )$ are versions of $\rho _{t_{0},t_{1}}(\mu )
$ and $\pi _{t_{0},t_{1}}(\mu )$.

Since $\rho _{t_{0},t_{1}}^{y_{(\cdot) }},$ $\pi _{t_{0},t_{1}}^{y_{(\cdot) }}$ can
be recast as time-inhomogenous dynamical systems, \emph{we can use them as
a basis for defining} $\rho _{t_{0},t_{1}}(\mu )$ and $\pi
_{t_{0},t_{1}}(\mu )$. As a result, it is indeed the case that $\rho
_{t_{0},t_{1}}(\mu )$ and $\pi _{t_{0},t_{1}}(\mu )$ can be viewed as RDSs.
Moreover, one can show that $\rho _{t_{0},\cdot }(\mu
)=\{\rho _{t_{0},t_{1}}(\mu ),\ t_{1}\geq t_{0}\}$~satisfies the evolution
equation 
\begin{equation}
\rho _{t_{0},t_{1}}(\mu )(\varphi )=\ \mu (\varphi
)+\int_{t_{0}}^{t_{1}}\rho _{t_{0},s}(\mu )(A\varphi )\,\mathrm{d}%
s+\int_{t_{0}}^{t_{1}}\rho _{t_{0},s}(\mu )(\varphi h^{\top })\mathrm{d}%
Y_{s},  \label{zakai}
\end{equation}%
for any $\varphi \in \mathcal{D}(A)$ and, in particular, that 
\begin{equation}
\rho _{t_{0},t_{1}}(\mu )(1)=1+\int_{t_{0}}^{t_{1}}\rho _{t_{0},s}(\mu
)(h^{\top })\mathrm{d}Y_{s}=1+\int_{t_{0}}^{t_{1}}\rho _{t_{0},s}(\mu
)\left( 1\right) \pi _{t_{0},s}(\mu )(h^{\top })\mathrm{d}Y_{s}.
\label{zakaimass}
\end{equation}%
From (\ref{zakai}) and (\ref{zakaimass}), one deduces that the ratio ${%
\rho _{t_{0},t_{1}}(\mu )(\varphi )/\rho _{t_{0},t_{1}}(\mu )(1)}$
satisfies (\ref{ks0}) with initial condition $\pi
_{t_{0},t_{0}}(\mu )=\mu $ and, therefore, by the uniqueness of the solution
of (\ref{ks0}), we obtain that the solution of (\ref%
{ks0}) can indeed be recast as an RDS. 

Finally, notice that $\pi _{0,t}(\pi _{0})=\pi_{t_{0},t}(\pi _{t_{0}})$ is indeed the FA process
we considerd throughout this paper. For
the particular case of $\mu =\pi _{t_{0}}$, formula (\ref{kallstrie})
is known as the Kallianpur-Striebel formula \cite{Moral.1996} and it is deduced directly
from the definition of the conditional expectation. For arbitrary $\mu $, (%
\ref{kallstrie}) serves as definition for the RDS, which
is then shown to be the solution of the evolution equation (\ref%
{ks0}) starting from $\mu $ at time $t_{0}$.
\begin{remark}
The map $\pi _{t_{0},t_{1}}^{y_{(\cdot) }}:\mathcal{P}(\mathbb{R}%
^{d})\rightarrow \mathcal{P}(\mathbb{R}^{d})$ is a continuous map when we
endow $\mathcal{P}(\mathbb{R}^{d})$ or, rather, the set of probability
measures with finite second moment), with the topology induced by the
Wasserstein metric.
\end{remark}


%
%

%

{\small 
\begin{acknowledgments}
It is a pleasure to thank Eviatar Bach for several useful discussions, for providing the final version of Figs.~\ref*{fig:Ghil1981} and \ref*{fig:FA_Cycle} and a draft of {\mg Appendix~I }, as well as for comments on the near-final manuscript. Alberto Carrassi and Pierre del Moral also read the near-final manuscript and made constructive suggestions.
M.G. acknowledges, however belatedly, the Nelder Fellowship of the Imperial College's Mathematics Department, \url{http://www.imperial.ac.uk/mathematics/research/opportunities/nelder-visiting-fellowships/fellows/professor-michael-ghil/} that put the two authors in closer contact in Spring 2014. Both authors are pleased to acknowledge the Institut Henri Poincar\'e's trimester on ``The Mathematics of Climate and the Environment” in Fall 2019, 
\url{http://www.ihp.fr/en/CEB/T3-2019},  which supported the real start of the collaboration leading to this paper. It is particularly gratifying to thank Franco Flandoli for an extended conversation during this trimester that set the two authors on the right path for the use of the tools in Crauel \& Flandoli\cite{Crauel.Flandoli.1994}. The comments of two reviewers --- one from the DA community, the other from the mathematical one --- have further improved the paper. D.C. was  partially supported by EU project STUOD - DLV-856408. The present paper is TiPES contribution \#  {xy}; this project has received funding from the EU Horizon 2020 research and innovation program under grant agreement No. 820970, and it helps support the work of M.G. Work on this paper has also been supported by the EIT Climate-KIC (grant no. 
190733); EIT Climate-KIC is supported by the European Institute of Innovation \& Technology (EIT), a body of the European Union.
\end{acknowledgments}
}

\bibliographystyle{abbrv}
\bibliography{DC+MG-v5b}

\begin{thebibliography}{100}

\bibitem{Abarbanel.ea.2017}
H.~D.~I. Abarbanel, S.~Shirman, D.~Breen, N.~Kadakia, D.~Rey, E.~Armstrong, and
  D.~Margoliash.
\newblock A unifying view of synchronization for data assimilation in complex
  nonlinear networks.
\newblock {\em Chaos: An Interdisciplinary Journal of Nonlinear Science},
  27(12):126802, 2017.

\bibitem{Arnold.1990}
L.~Arnold.
\newblock Stabilization by noise revisited.
\newblock {\em {ZAMM-Z. angew. Math. Mech.}}, 70(7):235--246, 1990.

\bibitem{Arnold.1998}
L.~Arnold.
\newblock {\em {Random Dynamical Systems}}.
\newblock Springer-Verlag, New York/Berlin, 1998.

\bibitem{Asch.Bocquet.2016}
M.~Asch, M.~Bocquet, and M.~Nodet.
\newblock {\em {Data Assimilation: Methods, Algorithms, and Applications}}.
\newblock SIAM, Philadelphia, PA, 2016.

\bibitem{atar2011}
R.~Atar.
\newblock Exponential decay rate of the filter's dependence on the initial
  distribution.
\newblock In {\em {The {O}xford Handbook of Nonlinear Filtering}}, pages
  299--318. Oxford Univ. Press, Oxford, 2011.

\bibitem{az1997}
R.~Atar and O.~Zeitouni.
\newblock Exponential stability for nonlinear filtering.
\newblock {\em Ann. Inst. H. Poincar\'{e} Probab. Statist.}, 33(6):697--725,
  1997.

\bibitem{Bach.Ghil.2022}
E.~Bach and M.~Ghil.
\newblock {A multi-model ensemble Kalman filter for data assimilation and
  forecasting}.
\newblock {\em arXiv E-print}, arXiv:2202.02272[physics, stat], 2022.

\bibitem{bc}
A.~Bain and D.~Crisan.
\newblock {\em {Fundamentals of Stochastic Filtering}}, volume~60 of {\em
  Stochastic Modelling and Applied Probability}.
\newblock Springer, New York, 2009.

\bibitem{bb2008}
D.~Bakry, F.~Barthe, P.~Cattiaux, and A.~Guillin.
\newblock A simple proof of the {P}oincar\'{e} inequality for a large class of
  probability measures including the log-concave case.
\newblock {\em Electron. Commun. Probab.}, 13:60--66, 2008.

\bibitem{Bakry.ea.2008}
D.~Bakry, P.~Cattiaux, and A.~Guillin.
\newblock Rate of convergence for ergodic continuous {M}arkov processes:
  {L}yapunov versus {P}oincar\'{e}.
\newblock {\em J. Funct. Anal.}, 254(3):727--759, 2008.

\bibitem{Bengtsson.ea.1981}
L.~Bengtsson, M.~Ghil, and E.~K{\"a}ll{\'e}n.
\newblock {\em {Dynamic Meteorology: Data Assimilation Methods}}.
\newblock Springer, 1981.

\bibitem{CrisanBeskosJasra}
A.~Beskos, D.~Crisan, and A.~Jasra.
\newblock On the stability of sequential {M}onte {C}arlo methods in high
  dimensions.
\newblock {\em Ann. Appl. Probab.}, 24(4):1396--1445, 2014.

\bibitem{bmm1}
A.~N. Bishop and P.~Del~Moral.
\newblock On the stability of matrix-valued {R}iccati diffusions.
\newblock {\em Electron. J. Probab.}, 24:Paper No. 84, 40, 2019.

\bibitem{bm2}
A.~N. Bishop and P.~Del~Moral.
\newblock On the stability of matrix-valued {R}iccati diffusions.
\newblock {\em Electron. J. Probab.}, 24:Paper No. 84, 40, 2019.

\bibitem{bmm2}
A.~N. Bishop and P.~Del~Moral.
\newblock An explicit {F}loquet-type representation of {R}iccati aperiodic
  exponential semigroups.
\newblock {\em Internat. J. Control}, 94(1):258--266, 2021.

\bibitem{bm}
A.~N. Bishop, P.~Del~Moral, K.~Kamatani, and B.~R\'{e}millard.
\newblock On one-dimensional {R}iccati diffusions.
\newblock {\em Ann. Appl. Probab.}, 29(2):1127--1187, 2019.

\bibitem{bm3}
A.~N. Bishop, P.~Del~Moral, and A.~Niclas.
\newblock A perturbation analysis of stochastic matrix {R}iccati diffusions.
\newblock {\em Ann. Inst. Henri Poincar\'{e} Probab. Stat.}, 56(2):884--916,
  2020.

\bibitem{Bjerknes.1904}
V.~Bjerknes.
\newblock {Das Problem der Wettervorhersage, betrachtet vom Standpunkte der
  Mechanik und der Physik}.
\newblock {\em Meteorologische Zeitschrift}, 21:1--7, 1904.

\bibitem{Bocquet.ea.2010}
M.~Bocquet, C.~A. Pires, and L.~Wu.
\newblock Beyond {Gaussian statistical modeling in geophysical data
  assimilation}.
\newblock {\em Monthly Weather Review}, 138(8):2997--3023, 2010.

\bibitem{Budhiraja2011}
A.~Budhiraja.
\newblock Feller and stability properties of the nonlinear filter.
\newblock In {\em The {Oxford Handbook of Nonlinear Filtering}}, pages
  352--373. Oxford Univ. Press, 2011.

\bibitem{Caraballo.Han.2017}
T.~Caraballo and X.~Han.
\newblock {\em {Applied Nonautonomous and Random Dynamical Systems: Applied
  Dynamical Systems}}.
\newblock Springer Science + Business Media, 2017.

\bibitem{Carrassi.ea.2008}
A.~Carrassi, M.~Ghil, A.~Trevisan, and F.~Uboldi.
\newblock Data assimilation as a nonlinear dynamical systems problem:
  {Stability and convergence of the prediction-assimilation system}.
\newblock {\em Chaos}, 18(2):023112, 2008.

\bibitem{Carrassietal2020}
A.~Carrassi, C.~Grudzien, M.~Bocquet, J.~Demaeyer, P.~Raanes, and S.~Vannitsem.
\newblock Data assimilation for chaotic systems.
\newblock In S.-K. Park and X.~Liang, editors, {\em {Data Assimilation for
  Atmospheric, Oceanic and Hydrological Applications}}. Springer Science \&
  Business Media, 2020.

\bibitem{CTDTU.2008}
A.~Carrassi, A.~Trevisan, L.~Descamps, O.~Talagrand, and F.~Uboldi.
\newblock {Controlling instabilities along a 3DVar analysis cycle by
  assimilating in the unstable subspace: a comparison with the EnKF}.
\newblock {\em Nonlinear Processes in Geophysics}, 15(4):503--521, 2008.

\bibitem{Carrassi.ea.2007}
A.~Carrassi, A.~Trevisan, and F.~Uboldi.
\newblock Adaptive observations and assimilation in the unstable subspace by
  breeding on the data-assimilation system.
\newblock {\em Tellus A: Dynamic Meteorology and Oceanography}, 59(1):101--113,
  2007.

\bibitem{Charney.ea.1969}
J.~G. Charney, M.~Halem, and R.~Jastrow.
\newblock Use of incomplete historical data to infer the present state of the
  atmosphere.
\newblock {\em {J. Atmos. Sci.}}, 26:1160--1163, 1969.

\bibitem{Charo.ea.2021}
G.~D. Char{\'o}, M.~D. Chekroun, D.~Sciamarella, and M.~Ghil.
\newblock Topological effects of noise on nonlinear dynamics.
\newblock {\em eprint arXiv:2010.09611v5 [nlin.CD]}, 2021.

\bibitem{Chavez.ea.2015}
M.~Chavez, M.~Ghil, and J.~Urrutia-Fucugauchi, editors.
\newblock {\em {Extreme Events: Observations, Modeling, and Economics}}, volume
  214.
\newblock John Wiley \& Sons, 2015.

\bibitem{Checkrounetal2011}
M.~D. Checkroun, E.~Simmonet, and M.~Ghil.
\newblock Stochastic climate dynamics: Random attractors and time-dependent
  invariant measures.
\newblock {\em Physica D: Nonlinear Phenomena}, 240:1685--1700, 2011.

\bibitem{cl2004}
P.~Chigansky and R.~Liptser.
\newblock Stability of nonlinear filters in nonmixing case.
\newblock {\em Ann. Appl. Probab.}, 14(4):2038--2056, 2004.

\bibitem{clh2011}
P.~Chigansky, R.~Liptser, and R.~Van~Handel.
\newblock Intrinsic methods in filter stability.
\newblock In {\em The {O}xford handbook of nonlinear filtering}, pages
  319--351. Oxford Univ. Press, Oxford, 2011.

\bibitem{Wei2}
C.~Cotter, D.~Crisan, D.~Holm, W.~Pan, and I.~Shevchenko.
\newblock Data assimilation for a quasi-geostrophic model with
  circulation-preserving stochastic transport noise.
\newblock {\em J. Stat. Phys.}, 179(5-6):1186--1221, 2020.

\bibitem{Wei3}
C.~Cotter, D.~Crisan, D.~D. Holm, W.~Pan, and I.~Shevchenko.
\newblock Numerically modeling stochastic {L}ie transport in fluid dynamics.
\newblock {\em {Multiscale Model. Simul.}}, 17(1):192--232, 2019.

\bibitem{Wei1}
C.~Cotter, D.~Crisan, D.~D. Holm, W.~Pan, and I.~Shevchenko.
\newblock A particle filter for stochastic advection by {L}ie transport: a case
  study for the damped and forced incompressible two-dimensional {E}uler
  equation.
\newblock {\em SIAM/ASA J. Uncertain. Quantif.}, 8(4):1446--1492, 2020.

\bibitem{Crauel.Flandoli.1994}
H.~Crauel and F.~Flandoli.
\newblock Attractors for random dynamical systems.
\newblock {\em Probability Theory and Related Fields}, 100(3):365--393, 1994.

\bibitem{Cressman.1959}
G.~P. Cressman.
\newblock An operational objective analysis system.
\newblock {\em Monthly Weather Review}, 87(10):367--374, 1959.

\bibitem{ch2008}
D.~Crisan and K.~Heine.
\newblock Stability of the discrete time filter in terms of the tails of noise
  distributions.
\newblock {\em J. Lond. Math. Soc. (2)}, 78(2):441--458, 2008.

\bibitem{uc1}
D.~Crisan and J.~M\'{\i}guez.
\newblock Uniform convergence over time of a nested particle filtering scheme
  for recursive parameter estimation in state-space {M}arkov models.
\newblock {\em Adv. Appl. Probab.}, 49(4):1170--1200, 2017.

\bibitem{CR}
D.~Crisan and B.~Rozovski\u{\i}, editors.
\newblock {\em {The {O}xford Handbook of Nonlinear Filtering}}.
\newblock Oxford University Press, Oxford, 2011.

\bibitem{Roisin.Beckers.2011}
B.~Cushman-Roisin and J.-M. Beckers.
\newblock {\em {Introduction to Geophysical Fluid Dynamics: Physical and
  Numerical Aspects}}.
\newblock Academic Press, 2nd Edition, 2011.
\newblock 875 pp.

\bibitem{Dalcher.Kalnay.1987}
A.~Dalcher and E.~Kalnay.
\newblock Error growth and predictability in operational {{ECMWF}} forecasts.
\newblock {\em Tellus A}, 39A(5):474--491, Oct. 1987.

\bibitem{Moral.1996}
P.~Del~Moral.
\newblock Nonlinear filtering using random particles.
\newblock {\em Theory of Probability \& Its Applications}, 40(4):690--701,
  1996.

\bibitem{uc2}
P.~Del~Moral.
\newblock {\em {Mean Field Simulation for {M}onte {C}arlo Integration}}, volume
  126 of {\em Monographs on Statistics and Applied Probability}.
\newblock CRC Press, Boca Raton, FL, 2013.

\bibitem{uc3}
P.~Del~Moral, A.~Doucet, and S.~S. Singh.
\newblock Uniform stability of a particle approximation of the optimal filter
  derivative.
\newblock {\em SIAM J. Control Optim.}, 53(3):1278--1304, 2015.

\bibitem{dmlm}
P.~Del~Moral and L.~Miclo.
\newblock On the stability of nonlinear {F}eynman-{K}ac semigroups.
\newblock {\em Ann. Fac. Sci. Toulouse Math. (6)}, 11(2):135--175, 2002.

\bibitem{mt1}
P.~Del~Moral and J.~Tugaut.
\newblock On the stability and the uniform propagation of chaos properties of
  ensemble {K}alman-{B}ucy filters.
\newblock {\em Ann. Appl. Probab.}, 28(2):790--850, 2018.

\bibitem{uc4}
P.~Del~Moral and J.~Tugaut.
\newblock Uniform propagation of chaos and creation of chaos for a class of
  nonlinear diffusions.
\newblock {\em Stoch. Anal. Appl.}, 37(6):909--935, 2019.

\bibitem{Dobrushin.1970}
R.~L. Dobrushin.
\newblock Prescribing a system of random variables by conditional
  distributions.
\newblock {\em Theory of Probability \& Its Applications}, 15(3):458--486,
  1970.

\bibitem{SMCM}
A.~Doucet, N.~de~Freitas, and N.~Gordon, editors.
\newblock {\em Sequential {M}onte {Carlo Methods in Practice}}.
\newblock Statistics for Engineering and Information Science. Springer-Verlag,
  New York, 2001.

\bibitem{Duane.ea.2017}
G.~S. Duane, C.~Grabow, F.~Selten, and M.~Ghil.
\newblock {Introduction to focus issue: Synchronization in large networks and
  continuous media{\textemdash}data, models, and supermodels}.
\newblock {\em Chaos}, 27(12):126601, dec 2017.

\bibitem{Duane.ea.2006}
G.~S. Duane, J.~J. Tribbia, and J.~B. Weiss.
\newblock Synchronicity in predictive modelling: a new view of data
  assimilation.
\newblock {\em {Nonlinear Processes in Geophysics}}, 13(6):601--612, 2006.

\bibitem{Farhat.ea.2015}
A.~Farhat, M.~S. Jolly, and E.~S. Titi.
\newblock {Continuous data assimilation for the 2D B{\'e}nard convection
  through velocity measurements alone}.
\newblock {\em Physica D: Nonlinear Phenomena}, 303:59--66, 2015.

\bibitem{Titi.ea.DA.2019}
A.~Farhat, E.~Lunasin, and E.~S. Titi.
\newblock {A data assimilation algorithm: The paradigm of the 3D Leray-$\alpha$
  model of turbulence}.
\newblock {\em Partial differential equations arising from physics and
  geometry}, pages 253--273, 2019.

\bibitem{flandoli3}
F.~Flandoli.
\newblock Some remarks on a statistical theory of turbulent flows.
\newblock In {\em Probabilistic Methods in Fluids}, pages 144--160. World Sci.
  Publ., River Edge, NJ, 2003.

\bibitem{Frank.Zhuk.2018}
J.~Frank and S.~Zhuk.
\newblock A detectability criterion and data assimilation for nonlinear
  differential equations.
\newblock {\em Nonlinearity}, 31(11):5235, 2018.

\bibitem{Gelb.1974}
A.~Gelb.
\newblock {\em {Applied Optimal Estimation}}.
\newblock M IT Press, Cambridge,M asachusetts, 1974.

\bibitem{Ghil.1997}
M.~Ghil.
\newblock Advances in sequential estimation for atmospheric and oceanic flows.
\newblock {\em J. Meteorol. Soc. Japan. Ser. II}, 75(1B):289--304, 1997.

\bibitem{Ghil.2015}
M.~Ghil.
\newblock A mathematical theory of climate sensitivity or, {How} to deal with
  both anthropogenic forcing and natural variability?
\newblock In C.-P. Chang, M.~Ghil, M.~Latif, and J.~Wallace, editors, {\em
  {Climate Change: Multidecadal and Beyond}}, volume~6, pages 31--52. World
  Scientific Publishing Co., Singapore, 2015.

\bibitem{GCS.2008}
M.~Ghil, M.~D. Chekroun, and E.~Simonnet.
\newblock Climate dynamics and fluid mechanics: natural variability and related
  uncertainties.
\newblock {\em Physica D: Nonlinear Phenomena}, 237(14--17):2111--2126, 2008.

\bibitem{Ghil.Chil.1987}
M.~Ghil and S.~Childress.
\newblock {\em {Topics in Geophysical Fluid Dynamics: Atmospheric Dynamics,
  Dynamo Theory, and Climate Dynamics}}.
\newblock Springer Science+Business Media, Berlin/Heidelberg, 1987.
\newblock {Reissued in pdf, 2012.}

\bibitem{Ghil.ea.1981}
M.~Ghil, S.~E. Cohn, J.~Tavantzis, K.~Bube, and E.~Isaacson.
\newblock Applications of estimation theory to numerical weather prediction.
\newblock In L.~Bengtsson, M.~Ghil, and E.~K\"all\'en, editors, {\em Dynamic
  Meteorology: Data Assimilation Methods}, pages 139--224. Springer, 1981.

\bibitem{Ghil.ea.1979}
M.~Ghil, M.~Halem, and R.~Atlas.
\newblock Time-continuous assimilation of remote-sounding data and its effect
  on weather forecasting.
\newblock {\em Monthly Weather Review}, 107:140--171, 1979.

\bibitem{Ghil.Luc.2020}
M.~Ghil and V.~Lucarini.
\newblock The physics of climate variability and climate change.
\newblock {\em Reviews of Modern Physics}, {in press}:{arXiv:1910.00583}, 2020.

\bibitem{Ghil.Mal.1991}
M.~Ghil and P.~Malanotte-Rizzoli.
\newblock Data assimilation in meteorology and oceanography.
\newblock {\em Adv. Geophys.}, 33({}):141--266, 1991.

\bibitem{Ghil.ea.1977}
M.~Ghil, B.~Shkoller, and V.~Yangarber.
\newblock A balanced diagnostic system compatible with a barotropic prognostic
  model.
\newblock {\em Monthly Weather Review}, 105(10):1223--1238, 1977.

\bibitem{Ghil.Todling.1996}
M.~Ghil and R.~Todling.
\newblock {Tracking atmospheric instabilities with the Kalman filter. Part II:
  Two-layer results}.
\newblock {\em Mon. Wea. Rev}, 124:2340--2352, 1996.

\bibitem{Ghil.ea.ExEv.2011}
M.~Ghil, P.~Yiou, S.~Hallegatte, B.~D. Malamud, P.~Naveau, A.~Soloviev,
  P.~Friederichs, V.~Keilis-Borok, D.~Kondrashov, V.~Kossobokov, O.~Mestre,
  C.~Nicolis, H.~W. Rust, P.~Shebalin, M.~Vrac, A.~Witt, and I.~Zaliapin.
\newblock Extreme events: dynamics, statistics and prediction.
\newblock {\em Nonlinear Processes in Geophysics}, 18(3):295--350, 2011.

\bibitem{Gottwald.Reich.2021}
G.~A. Gottwald and S.~Reich.
\newblock Combining machine learning and data assimilation to forecast
  dynamical systems from noisy partial observations.
\newblock {\em Chaos: An Interdisciplinary Journal of Nonlinear Science},
  31(10):101103, oct 2021.

\bibitem{Halem.ea.1982}
M.~Halem, E.~Kalnay, W.~E. Baker, and R.~Atlas.
\newblock An assessment of the {FGGE satellite observing system during SOP-1}.
\newblock {\em Bulletin of the American Meteorological Society},
  63(4):407--426, 1982.

\bibitem{uc5}
K.~Heine and D.~Crisan.
\newblock Uniform approximations of discrete-time filters.
\newblock {\em Adv. in Appl. Probab.}, 40(4):979--1001, 2008.

\bibitem{Hoteit.ea.2015}
I.~Hoteit, D.-T. Pham, M.~E. Gharamti, and X.~Luo.
\newblock Mitigating observation perturbation sampling errors in the stochastic
  {EnKF}.
\newblock {\em Monthly Weather Review}, 143(7):2918--2936, 2015.

\bibitem{Houtekamer.2016}
P.~Houtekamer and F.~Zhang.
\newblock {Review of the ensemble Kalman filter for atmospheric data
  assimilation}.
\newblock {\em Monthly Weather Review}, 144(12):4489--4532, 2016.

\bibitem{Ide.ea.1997}
K.~Ide, P.~Courtier, M.~Ghil, and A.~C. Lorenc.
\newblock {Unified notation for data assimilation: Operational, sequential and
  variational, in {\it Special Issue on Data Assimilation in Meteology and
  Oceanography: Theory and Practice}}.
\newblock {\em Journal of the Meteorological Society of Japan. Ser. II},
  75(1B):181--189, 1997.

\bibitem{ikedawatanabe}
N.~Ikeda and S.~Watanabe.
\newblock {\em {Stochastic Differential Equations and Diffusion Processes}},
  volume~24 of {\em North-Holland Mathematical Library}.
\newblock North-Holland Publishing Co., Amsterdam-New York; Kodansha, Ltd.,
  Tokyo, 1981.

\bibitem{Jazwinski.1970}
A.~H. Jazwinski.
\newblock {\em {Stochastic Processes and Filtering Theory}}.
\newblock Dover Publ., {1970, reprinted by Courier Corp. in 2007}.

\bibitem{Kalman.1960}
R.~Kalman.
\newblock A new approacht to linear filtering and prediction problems.
\newblock {\em ASME J. Basic Eng.}, 82D:35--45, 1960.

\bibitem{Kalman.Bucy.1961}
R.~Kalman and R.~Bucy.
\newblock New results in linear filtering and prediction theory.
\newblock {\em ASME J. Basic Eng.}, 83D:95--108, 1961.

\bibitem{Kalnay.2003}
E.~Kalnay.
\newblock {\em Atmospheric {Modeling, Data Assimilation and Predictability}}.
\newblock Cambridge University Press, Cambridge, UK, 2003.

\bibitem{KantasBeskosJasra}
N.~Kantas, A.~Beskos, and A.~Jasra.
\newblock Sequential {M}onte {C}arlo methods for high-dimensional inverse
  problems: a case study for the {N}avier-{S}tokes equations.
\newblock {\em SIAM/ASA J. Uncertain. Quantif.}, 2(1):464--489, 2014.

\bibitem{Kantorovich.1942}
L.~V. Kantorovich.
\newblock On the translocation of masses.
\newblock {\em Journal of Mathematical Sciences}, 133(4):1381--1382, 2006.
\newblock {originally published in Doklady Akademii Nauk SSSR, 37 (7–8),
  199--201 (1942).}

\bibitem{KaratzasShreve}
I.~Karatzas and S.~E. Shreve.
\newblock {\em {Brownian Motion and Stochastic Calculus}}, volume 113 of {\em
  Graduate Texts in Mathematics}.
\newblock Springer-Verlag, New York, 1988.

\bibitem{km}
T.~Karvonen, S.~Bonnabel, E.~Moulines, and S.~S\"{a}rkk\"{a}.
\newblock On stability of a class of filters for nonlinear stochastic systems.
\newblock {\em SIAM J. Control Optim.}, 58(4):2023--2049, 2020.

\bibitem{kun}
H.~Kunita.
\newblock Asymptotic behavior of the nonlinear filtering errors of {M}arkov
  processes.
\newblock {\em J. Multivariate Anal.}, 1:365--393, 1971.

\bibitem{Kunita.1986}
H.~Kunita.
\newblock Stochastic flows and stochastic partial differential equations.
\newblock In {\em Proceedings of the {I}nternational {C}ongress of
  {M}athematicians, {V}ol. 1, 2 ({B}erkeley, {C}alif., 1986)}, pages
  1021--1031. Amer. Math. Soc., Providence, RI, 1987.

\bibitem{ks}
H.~Kwakernaak and R.~Sivan.
\newblock {\em {Linear Optimal Control Systems}}.
\newblock Wiley-Interscience [John Wiley \& Sons], New York-London-Sydney,
  1972.

\bibitem{Stuart.ea.2015}
K.~Law, A.~Stuart, and K.~Zygalakis.
\newblock {\em {Data Assimilation: A Mathematical Introduction}}.
\newblock {Springer, Cham, Switzerland}, 2015.

\bibitem{Lawson.Hansen.2004}
W.~Lawson and J.~Hansen.
\newblock Implications of stochastic and deterministic filters as
  ensemble-based data assimilation methods in varying regimes of error growth.
\newblock {\em Monthly Weather Review}, 132(8):1966--1981, 2004.

\bibitem{uc6}
F.~Le~Gland and N.~Oudjane.
\newblock Stability and uniform approximation of nonlinear filters using the
  {H}ilbert metric and application to particle filters.
\newblock {\em Ann. Appl. Probab.}, 14(1):144--187, 2004.

\bibitem{Leeuwen.2020}
P.~J. Leeuwen.
\newblock {A consistent interpretation of the stochastic version of the
  Ensemble Kalman Filter}.
\newblock {\em Quarterly Journal of the Royal Meteorological Society},
  146(731):2815--2825, jun 2020.

\bibitem{Leith.1974}
C.~E. Leith.
\newblock Theoretical skill of {{Monte Carlo forecasts}}.
\newblock {\em Monthly Weather Review}, 102(6):409--418, June 1974.

\bibitem{Leith.1978}
C.~E. Leith.
\newblock Objective methods for weather prediction.
\newblock {\em Annual Review of Fluid Mechanics}, 10(1):107--128, 1978.

\bibitem{lr}
W.~Liu and M.~R\"{o}ckner.
\newblock {\em Stochastic partial differential equations: an introduction}.
\newblock Universitext. Springer, Cham, 2015.

\bibitem{Lorenz.1963a}
E.~N. Lorenz.
\newblock Deterministic nonperiodic flow.
\newblock {\em Journal of the Atmospheric Sciences}, 20:130--141, 1963.

\bibitem{Lorenz.Book.1967}
E.~N. Lorenz.
\newblock {\em {The Nature and Theory of the General Circulation of the
  Atmosphere}}, volume 218.
\newblock World Meteorological Organization Geneva, 1967.

\bibitem{Lorenz.1982}
E.~N. Lorenz.
\newblock Atmospheric predictability experiments with a large numerical model.
\newblock {\em Tellus}, 34(6):505--513, 1982.

\bibitem{Lorenz.1984}
E.~N. Lorenz.
\newblock {Irregularity: A fundamental property of the atmosphere}.
\newblock {\em Tellus A}, 36(2):98--110, 1984.

\bibitem{Malguzzi.ea.1990}
P.~Malguzzi, A.~Trevisan, and A.~Speranza.
\newblock Statistics and predictability for an intermediate dimensionality
  model of the baroclinic jet.
\newblock {\em {Annales Geophysicae. Atmospheres, Hydrospheres and Space
  Sciences}}, 8(1):29--35, 1990.

\bibitem{Monge.1781}
G.~Monge.
\newblock {Mémoire sur la théorie des déblais et des remblais}.
\newblock {\em Histoire de l’Académie Royale des Sciences}, pages 666--704,
  1781.

\bibitem{Nicolis.ea.2009}
C.~Nicolis, R.~A.~P. Perdigao, and S.~Vannitsem.
\newblock Dynamics of prediction errors under the combined effect of initial
  condition and model errors.
\newblock {\em Journal of the Atmospheric Sciences}, 66(3):766--778, 2009.

\bibitem{nualart}
D.~Nualart.
\newblock {\em The {M}alliavin calculus and related topics}.
\newblock Probability and its Applications (New York). Springer-Verlag, Berlin,
  second edition, 2006.

\bibitem{op}
D.~Ocone and E.~Pardoux.
\newblock Asymptotic stability of the optimal filter with respect to its
  initial condition.
\newblock {\em SIAM J. Control Optim.}, 34(1):226--243, 1996.

\bibitem{Panaretos.Zemel.2020}
V.~M. Panaretos and Y.~Zemel.
\newblock {\em {An Invitation to Statistics in Wasserstein Space}}.
\newblock Springer Nature, 2020.

\bibitem{Panofsky.1949}
R.~A. Panofsky.
\newblock Objective weather-map analysis.
\newblock {\em Journal of Meteorology}, 6(6):386--392, 1949.

\bibitem{picard1991}
J.~Picard.
\newblock Efficiency of the extended {K}alman filter for nonlinear systems with
  small noise.
\newblock {\em {SIAM J. Appl. Math.}}, 51(3):843--885, 1991.

\bibitem{Potthast}
R.~Potthast, A.~Walter, and A.~Rhodin.
\newblock A localized adaptive particle filter within an operational {NWP
  framework}.
\newblock {\em Monthly Weather Review}, 147:345--362, 2019.

\bibitem{rc}
S.~Reich and C.~Cotter.
\newblock {\em Probabilistic {Forecasting and {B}ayesian Data Assimilation}}.
\newblock Cambridge University Press, New York, 2015.

\bibitem{rt1996}
G.~O. Roberts and R.~L. Tweedie.
\newblock Exponential convergence of {L}angevin distributions and their
  discrete approximations.
\newblock {\em Bernoulli}, 2(4):341--363, 1996.

\bibitem{Robin.ea.2017}
Y.~Robin, P.~Yiou, and P.~Naveau.
\newblock {Detecting changes in forced climate attractors with Wasserstein
  distance}.
\newblock {\em Nonlinear Processes in Geophysics}, 24:393--405, 2017.

\bibitem{rl}
B.~L. Rozovsky and S.~V. Lototsky.
\newblock {\em Stochastic evolution systems}, volume~89 of {\em Probability
  Theory and Stochastic Modelling}.
\newblock Springer, Cham, 2018.
\newblock Linear theory and applications to non-linear filtering, Second
  edition of [ MR1135324].

\bibitem{Savijarvi.1995}
H.~Savij\"arvi.
\newblock Error growth in a large numerical forecast system.
\newblock {\em Monthly Weather Review}, 123(1):212--221, jan 1995.

\bibitem{Simmons.ea.1995}
A.~J. Simmons, R.~Mureau, and T.~Petroliagis.
\newblock Error growth and estimates of predictability from the {ECMWF}
  forecasting system.
\newblock {\em Quarterly Journal of the Royal Meteorological Society},
  121(527):1739--1771, oct 1995.

\bibitem{stannat2005}
W.~Stannat.
\newblock Stability of the filter equation for a time-dependent signal on
  {$\Bbb R^d$}.
\newblock {\em Appl. Math. Optim.}, 52(1):39--71, 2005.

\bibitem{stannat2011}
W.~Stannat.
\newblock Stability of the optimal filter for nonergodic signals---a
  variational approach.
\newblock In {\em The {O}xford handbook of nonlinear filtering}, pages
  374--399. Oxford Univ. Press, Oxford, 2011.

\bibitem{Royer.1993}
R.~Stroe and J.~F. Royer.
\newblock Comparison of different error growth formulas and predictability
  estimation in numerical extended-range forecasts.
\newblock {\em Annales Geophysicae}, 11(4):296--316, 1993.

\bibitem{Trevisan.ea.2010}
A.~Trevisan, M.~D'Isidoro, and O.~Talagrand.
\newblock Four-dimensional variational assimilation in the unstable subspace
  and the optimal subspace dimension.
\newblock {\em Quarterly Journal of the Royal Meteorological Society},
  136(647):487--496, 2010.

\bibitem{Trevisan.ea.1992}
A.~Trevisan, P.~Malguzzi, and M.~Fantini.
\newblock {On Lorenz's law for the growth of large and small errors in the
  atmosphere}.
\newblock {\em Journal of the Atmospheric Sciences}, 49(8):713--719, apr 1992.

\bibitem{Trevisan.Uboldi.2004}
A.~Trevisan and F.~Uboldi.
\newblock Assimilation of standard and targeted observations within the
  unstable subspace of the observation--analysis--forecast cycle system.
\newblock {\em Journal of the Atmospheric Sciences}, 61(1):103--113, 2004.

\bibitem{Handel.2009}
R.~van Handel.
\newblock {Uniform observability of hidden Markov models and filter stability
  for unstable signals}.
\newblock {\em The Annals of Applied Probability}, 19(3), 2009.

\bibitem{Leeuwen.2009}
P.~J. Van~Leeuwen.
\newblock Particle filtering in geophysical systems.
\newblock {\em Monthly Weather Review}, 137(12):4089--4114, 2009.

\bibitem{Leeuwen2010}
P.~J. van Leeuwen.
\newblock Nonlinear data assimilation in geosciences: an extremely efficient
  particle filter.
\newblock {\em Quarterly Journal of the Royal Meteorological Society}, 136,
  2010.

\bibitem{Leeuwen.ea.2015}
P.~J. Van~Leeuwen, Y.~Cheng, and S.~Reich.
\newblock {\em {Nonlinear Data Assimilation}}.
\newblock Springer, 2015.

\bibitem{VetraCarvalho-Leeuwen2018}
S.~Vetra-Carvalho, P.~J. van Leeuwen, L.~Nerger, A.~Barth, M.~U. Altaf,
  P.~Brasseur, P.~Kirchgessner, and J.~M. Beckers.
\newblock State-of-the-art stochastic data assimilation methods for
  high-dimensional {non-Gaussian problems}.
\newblock {\em Tellus A}, 70(1):1--43, 2018.

\bibitem{Villani.2009}
C.~Villani.
\newblock {\em {Optimal Transport: Old and New}}.
\newblock Springer-Verlag, Berlin Heidelberg, Germany, 2009.

\bibitem{Vissio.ea.2020}
G.~Vissio, V.~Lembo, V.~Lucarini, and M.~Ghil.
\newblock Evaluating the performance of climate models based on {Wasserstein
  distance}.
\newblock {\em Geophysical Research Letters}, 47(21), oct 2020.

\bibitem{Wasserstein.1969}
L.~Wasserstein.
\newblock Markov processes with countable state space describing large systems
  of automata.
\newblock {\em Problemy Peredachi Informatsii}, 5:64--73, {1969 (in Russian)}.

\bibitem{Wax.1954}
N.~Wax.
\newblock {\em {Selected Papers on Noise and Stochastic Processes}}.
\newblock Courier Dover Publications, 1954.

\end{thebibliography}

\end{document}